\documentclass[12pt]{article}

\usepackage[utf8]{inputenc}
\usepackage[T1]{fontenc}

\usepackage{amsmath,mleftright,mathtools,xparse,
            amssymb,amsthm,nicefrac,xcolor,
            enumerate,color,comment,enumitem,
            mathrsfs,bbm,stmaryrd,geometry,etoolbox
            }
            
\mleftright   

\usepackage[sort,nocompress]{cite}

\usepackage[colorinlistoftodos]{todonotes}
\tikzset{notestyleraw/.append style={align=justify}}
                        
\usepackage[english]{datetime2}
            
\usepackage[colorlinks=true]{hyperref}
            
\newcommand{\N}{\ensuremath{\mathbb{N}}}
\newcommand{\Z}{\ensuremath{\mathbb{Z}}}

\newcommand{\R}{\ensuremath{\mathbb{R}}}
\newcommand{\C}{\ensuremath{\mathbb{C}}}
\newcommand{\1}{\mathbbm{1}}
\newcommand{\sSpace}{\mathbb{S} \cfadd{def:seq}}
\newcommand{\semigroup}{S \cfadd{def:semigroup}}
\newcommand{\id}{\mathbb{I} \cfadd{def:id}}
\newcommand{\Binom}[2]{{\textstyle\binom{#1}{#2}}}
\newcommand{\ellTwo}{\ell^2(\Z) \cfadd{def:l2}}
\newcommand{\lap}{\Delta \cfadd{def:lap}}
\newcommand{\gamfn}{\Gamma \cfadd{def:gamma}}
\newcommand{\preFrac}{(-\Delta)^s \cfadd{def:frac_lap_pre}}
\newcommand{\fracLap}{(-\Delta)^s \cfadd{def:frac_lap}}
\newcommand{\Kern}{K \cfadd{def:kern}}
\newcommand{\real}{\mathfrak{R}}
\newcommand{\logNorm}[1]{\mu(#1) \cfadd{def:log}}
\newcommand{\dini}[1]{D^+_{#1} \cfadd{def:dini}}

\newcommand{\SmallSum}{\textstyle\sum\limits}

\newcommand{\induct}{\dashrightarrow}
\newcommand{\with}{\curvearrowleft}
\newcommand{\lrSpace}{\ensuremath{\mkern-1.5mu}}

\DeclarePairedDelimiter{\pr}{(}{)}
\DeclarePairedDelimiter{\br}{[}{]}
\DeclarePairedDelimiter{\cu}{\{}{\}}
\DeclarePairedDelimiter{\abs}{\lvert}{\rvert}
\DeclarePairedDelimiter{\norm}{\lVert}{\rVert}
\DeclarePairedDelimiter{\floor}{\lfloor}{\rfloor}
\DeclarePairedDelimiter{\vt}{\langle}{\rangle}
\DeclarePairedDelimiter{\inner}{\langle}{\rangle \cfadd{def:inner}}
\DeclarePairedDelimiter{\semi}{[}{] \cfadd{def:semi}}

\DeclarePairedDelimiter{\normm}{\lVert}{\rVert_2 \cfadd{def:l2}}

\allowdisplaybreaks

\usepackage[sort,capitalize]{cleveref}

\crefformat{equation}{(#2#1#3)}
\crefname{equation}{}{}
\crefname{enumi}{item}{items}
\crefname{subsection}{Subsection}{Subsections}

\newtheorem{theorem}{Theorem}[section]

\newtheorem{corollary}[theorem]{Corollary}

\theoremstyle{definition}
\newtheorem{definition}[theorem]{Definition}

\numberwithin{equation}{section}

\ExplSyntaxOn

\seq_new:N \g_cflist_loaded
\seq_new:N \g_cflist_pending

\NewDocumentCommand{\cfadd}{ m }
{
  \seq_if_in:NnF \g_cflist_loaded { #1 } {
    \seq_if_in:NnF \g_cflist_pending { #1 } {
      \seq_gput_right:Nn \g_cflist_pending { #1 }
    }
  }
}

\NewDocumentCommand{\cfconsiderloaded}{ m }{
  \seq_gput_right:Nn \g_cflist_loaded {#1}
}

\NewDocumentCommand{\cfremove}{ m }
{
  \seq_gremove_all:Nn \g_cflist_pending { #1 }
}

\NewDocumentCommand{\cfload}{ o }
{
  \seq_if_empty:NTF \g_cflist_pending {\unskip} {
    (cf.\ \cref{\seq_use:Nn \g_cflist_pending {,}})\IfValueTF{#1}{#1~}{\unskip}
    \seq_gconcat:NNN \g_cflist_loaded \g_cflist_loaded \g_cflist_pending
    \seq_gclear:N \g_cflist_pending
  }
}

\NewDocumentCommand{\cfclear} {} {
  \seq_gclear:N \g_cflist_loaded
  \seq_gclear:N \g_cflist_pending
}

\NewDocumentCommand{\cfout}{ o }
{
  \seq_if_empty:NTF \g_cflist_pending {\unskip} {
    (cf.\ \cref{\seq_use:Nn \g_cflist_pending {,}})\IfValueTF{#1}{#1~}{\unskip}
    \seq_gclear:N \g_cflist_pending
  }
}

\NewDocumentCommand{\ifnocf} { m } {
  \seq_if_empty:NT \g_cflist_pending { #1 }
}

\NewDocumentEnvironment {athm} {m m} {%
\begin{#1}\label{#2}\global\def\loc{#2}%
}{%
\end{#1}%
}

\NewDocumentEnvironment{aproof} {} {%
\begin{proof}[Proof~of~\cref{\loc}]%
}{%
\finishproofthus
\end{proof}%
}

\ExplSyntaxOff

\newcommand{\finishproofthus}{The proof of \cref{\loc} is thus complete.}

\begin{document}

\title{A series representation of the discrete fractional Laplace operator of arbitrary order}

\author{
Tiffany Frug{\'e} Jones$^1$,
Evdokiya Georgieva Kostadinova$^2$,\\
Joshua Lee Padgett$^{3,4}$, and Qin Sheng$^{4,5}$
\bigskip
\\
\small{$^1$ Department of Mathematics, University of Arizona,}
\vspace{-0.1cm}\\
\small{Tucson, Arizona 85721, USA, e-mail: \texttt{tnjones@math.arizona.edu}}
\smallskip
\\
\small{$^2$ Center for Astrophysics, Space Physics, and Engineering Research,}
\vspace{-0.1cm}\\
\small{Baylor University, Waco, Texas 76798, USA, e-mail: \texttt{eva\_kostadinova@baylor.edu}}
\smallskip
\\
\small{$^3$ Department of Mathematical Sciences, University of Arkansas,}
\vspace{-0.1cm}\\
\small{Fayetteville, Arkansas 72701, USA, e-mail: \texttt{padgett@uark.edu}}
\smallskip
\\
\small{$^4$ Center for Astrophysics, Space Physics, and Engineering Research,}
\vspace{-0.1cm}\\
\small{Baylor University, Waco, Texas 76798, USA}
\smallskip
\\
\small{$^5$ Department of Mathematics, Baylor University,}
\vspace{-0.1cm}\\
\small{Waco, Texas, USA, e-mail: \texttt{qin\_sheng@baylor.edu}}
\smallskip
}

\date{\today}

\maketitle

\begin{abstract}
Although fractional powers of non-negative operators have received much attention in recent years, 
there is still little known about their behavior if real-valued exponents are greater than one. 
In this article, we define and study the discrete fractional Laplace operator of arbitrary real-valued positive order. 
A series representation of the discrete fractional Laplace operator
for positive non-integer powers is developed. Its convergence to a series representation of
a known case of positive integer powers is proven as the power tends to the integer value. 
Furthermore, 
we show that the new representation for arbitrary real-valued positive powers of the discrete Laplace 
operator is consistent with existing theoretical results.
\end{abstract}

\tableofcontents


\section{Introduction}
\label{sec:intro}

Due to its wide array of applications in multi-physical sciences, the construction and approximation of 
fractional powers of the Laplace operator have been of great interest for nearly a century
(cf., e.g., \cite{pozrikidis2018fractional,%
bucur2016nonlocal,meerschaert2011stochastic,%
vazquez2017mathematical,lischke2020fractional} and references therein).
Conventionally, only powers of the order $s\in(0,1)$ are considered, and in this case, 
one may define the fractional Laplace operator applied to a smooth enough function in a natural way.
Specifically, for $d \in \N = \{1,2,3,\ldots\}$, $s\in(0,1)$ let $u \colon \R^d \to \R$ be a smooth function, and for every 
$\varepsilon\in(0,\infty)$, $x\in\R^d$ let $B_\varepsilon(x)$ be the $d$-dimensional ball of radius 
$\varepsilon$ centered at $x$ (with respect to the typical topology of $\R^d$). Then for every $x\in\R^d$ 
we can define the $s$-order fractional Laplace operator applied to $u$ at $x$ as
\begin{equation}\label{eq:caff}
\pr[\big]{ (-\Delta)^s u }(x) = 
c_{d,s} \lim_{\varepsilon\to 0^+} \left[ \int_{\R^d \backslash B_\varepsilon(x) } 
\frac{ u(x) - u(y) }{ \abs{ x - y }^{d+2s} } \, dy \right],
\end{equation}
where $c_{d,s} \in [0,\infty)$ is a known normalization constant.

It is worth noting that the recent rapid increase in interest in the fractional Laplace operator 
is also due to the seminal work of Caffarelli and Silvestre \cite{caffarelli2007extension}.
In their work, it was shown that one may study the non-local operator given by \cref{eq:caff} 
via the Dirichlet-to-Neumann operator associated with a particular extension problem posed in 
$\R^d \times [0,\infty)$ 
(albeit, one trades the non-locality for a problem posed in a higher dimension which is either singular or 
degenerate depending upon the value of $s\in(0,1)$).
The employed Dirichlet-to-Neumann operator is a particular example of the Poincar{\'e}-Stecklov operator (cf., e.g., \cite{poin}).
For a fixed domain, the Poincar{\'e}-Stecklov operator is known to map the boundary values of a harmonic function to the normal derivative values of the same harmonic function on the same boundary.
We can summarize the results of Caffarelli and Silvestre (cf., e.g., \cite[Eq.\ (3.1)]{caffarelli2007extension}) as follows.
Let $d \in \N$, $s\in(0,1)$, let $u \colon \R^d \to \R$ be a smooth function, and let $v \colon [0,\infty) \times \R^d \to \R$ satisfy for all $x\in\R^d$ that $v(0,x) = u(x)$ and for all $t\in(0,\infty)$, $x\in\R^d$ that
\begin{equation}\label{eq:caff1}
\left(\tfrac{\partial^2}{\partial t^2} v\right)(t,x) + \tfrac{1-2s}{t} \left(\tfrac{\partial}{\partial t} v\right)(t,x) + (\Delta_x v)(t,x) = 0.
\end{equation}
Then there exists $c \in [0,\infty)$ such that for all $x\in\R^d$ it holds that
\begin{equation}\label{eq:caff2}
\pr[\big]{ (-\Delta)^s u }(x) = c \left[ \lim_{t\to 0^+} t^{1-2s} \left(\tfrac{\partial}{\partial t} v\right)(t,x) \right] .
\end{equation}
Interestingly, the constant $c\in [0,\infty)$ in \cref{eq:caff2} depends only upon the parameter $s \in (0,1)$ and not upon $d\in\N$.
More importantly, this demonstrates that one may trade out the highly non-local problem given by \cref{eq:caff} 
for the local problem given by \cref{eq:caff1,eq:caff2}.
This technique has also been recently further generalized to cases of arbitrary non-negative operators defined on 
Banach spaces \cite{padgett2020analysis,meichsner2017fractional,%
gale2013extension,arendt2018fractional,%
meichsner2020harmonic}.

While the above formulations (i.e., \cref{eq:caff,eq:caff1,eq:caff2}) may be utilized to provide insights into 
a continuous fractional Laplace operator with order $s\in(0,1)$, they cannot be generalized to provide any
insight into the discrete case or the case where $s \in (0,\infty)$.
The discrete case is a natural consideration as it arises in the study of numerous physically relevant 
phenomena (cf., e.g., \cite{PhysRevResearch.2.043375,%
padgett2020anomalous,kostadinova2017delocalization}
and references therein)
and also in an attempt to numerically approximate \cref{eq:caff}.
The consideration of a truly discrete case---that is, the case which is the fractional power of the discrete 
Laplace operator rather than a direct approximation of \cref{eq:caff}---was originally studied by 
Ciaurri et al.\ \cite{ciaurri2018nonlocal}.
By employing the basic language of semigroups (e.g., a special case of  
Ciaurri et al.\ \cite[Eq.\ (1)]{ciaurri2017harmonic} combined with, e.g.,
Padgett \cite[Theorem 2.1]{padgett2020analysis})
Ciaurri et al.\ were able to develop the first series representation for the discrete fractional 
Laplace operator of order $s\in(0,1)$ (cf.\ \cref{def:frac_lap_pre}, for clarity).
Moreover, it was shown that this formulation did converge to the continuous case via 
adaptive mesh refinements (cf.\ Ciaurri et al.\ \cite[Theorems 1.7 and 1.8]{ciaurri2017harmonic}).
However, it important to note that while this aforementioned convergence was observed, it is the case that the series representation developed by Ciaurri et al.\ is an \emph{exact} representation and not a numerical approximation.

The consideration of higher-order fractional Laplace operators has recently received increased attention in continuous cases
(cf., e.g., \cite{chen2018extension,yang2013higher,%
ros2014local,garcia2019strong,felli2018unique}).
But to the authors' knowledge, the only study in discrete cases has been carried out by Padgett et al.\
\cite{padgett2020anomalous}.
Rectifying this aforementioned gap in theory is the primary goal of this article (although the applicability of such 
derivations in the study of localization will be outlined in \cref{sec:motivate} below). 
To this end, a series representation of the discrete fractional Laplace operator 
of order $s\in(0,\infty)$ is implemented.
This development is illustrated through \cref{th:main}, which is also a partial description of the main result of this 
article focused on the case of positive non-integer powers of the discrete Laplace operator.

\begin{athm}{theorem}{th:main}
Let $m \in \N$, $s\in(m-1,m)$, 
let $\Z = \{\ldots,-2,-1,0,1,2,\ldots\}$, 
let $\R$ be the real number field, let $\ellTwo$ be the set of all $w \colon \Z \to \R$ which satisfy that $\sum_{k\in\Z} \abs{w(k)}^2 < \infty$, 
let $-\Delta \colon \ellTwo \to \ellTwo$ satisfy for all $w \in \ellTwo$, $n\in\Z$ that $(-\Delta w)(n) = 2w(n) - w(n-1) - w(n+1)$, 
let $u \in \ellTwo$,
and
let\footnote{Note that we define integer powers of $-\lap \colon \ellTwo \to \ellTwo$ inductively. That is, we have for all $k\in\N_0 = \N \cup \{0\}$, $w \in \ellTwo$, $n\in\Z$ that if $k=0$ it holds that $( (-\lap)^k w)(n) = -w(n) $ and if $k \in \N$ it holds that $( (-\lap)^k w)(n) = ( -\lap (-\lap)^{k-1} w)(n)$.}
$v \colon [0,\infty) \times \Z \to \R$ satisfy for all $t\in(0,\infty)$, $n\in\Z$ that $v(0,n) = ( (-\Delta)^{m-1} u)(n)$ and
\begin{equation}\label{ex:extend}
\left(\tfrac{\partial^2}{\partial t^2} v\right)(t,x) + \tfrac{1-2(s-m+1)}{t} \left(\tfrac{\partial}{\partial t} v\right)(t,x) + \left(\Delta v\right)(t,x) = 0.
\end{equation}
Then
\begin{enumerate}[label=(\roman*)]
\item\label{th:main_i1} there exists $c \in [0,\infty)$ such that for all $n\in\Z$ it holds that
\begin{equation}\label{eq:1_5}
\pr[\big]{ (-\Delta)^{s-m+1} (-\Delta)^{m-1} u }(n) =
\pr[\big]{ (-\Delta)^s u }(n) = c \left[ \lim_{t\to 0^+} t^{1-2(s-m+1)} \left(\tfrac{\partial}{\partial t} v\right)(t,n) \right]
\end{equation}
and
\item\label{th:main_i2} there exists $K \colon \Z \to \R$, $C \in [0,\infty)$ such that for all $n\in\Z\backslash\{0\}$ it holds that
$\abs{K(n)} \le C \abs{n}^{-(1+2s)}$ and for all $n\in\Z$ it holds that $K(-n) = K(n)$ and
\begin{equation}
\pr[\big]{ (-\Delta)^s u }(n) = \sum_{k\in\Z} K(k) \pr[\big]{ u(n) - u(n-k) } .
\end{equation}
\end{enumerate}
\end{athm}

We now provide some clarifying remarks regarding the objects in \cref{th:main}.
In \cref{th:main} we intend to construct an exact series representation of the so-called co-normal derivative of the function $v(0,\cdot) \colon \Z \to \R$.
The positive real number $s\in(0,\infty)$ describes the fractional power of the discrete Laplace operator, the positive integer $m$ describes the smallest positive integer that is greater than or equal to $s \in (0,\infty)$, and the set $\ellTwo$ is the standard Hilbert space of square-summable sequences defined on the integers.
The operator $-\Delta \colon \ellTwo \to \ellTwo$ is the standard one-dimensional discrete Laplace operator and is the primary object used in the construction of the desired series representation.
The function $v \colon [0,\infty) \times \Z \to \R$ is the solution to the extension problem in \cref{ex:extend} and the trace of this function coincides with the given square-summable function $u \colon \Z \to \R$ to which we are applying the discrete fractional Laplace operator of order $s \in (0,\infty)$.

Let us now provide some clarifying remarks regarding the results in \cref{th:main}.
\Cref{th:main_i1} of \cref{th:main} above is a direct consequence of combining \cref{def:frac_lap} and Padgett \cite[Theorem 2.1]{padgett2020analysis} (applied for every $n\in\Z$ with $s \with s-m+1$, $A \with \Delta$, $u_0 \with ( (-\Delta)^{m-1} u)(n)$, $(u(t))_{t\in[0,\infty)} \with (v(t,n))_{t\in[0,\infty)}$ in the notation of Padgett \cite[Theorem 2.1]{padgett2020analysis}).
See the beginning of \cref{sec:2} for an explanation of this ``applied with'' notation (i.e., the symbol ``$\with$'').
The right-hand side of \cref{eq:1_5} is not considered in detail, herein, as it is an elementary consequence of, e.g., Padgett \cite[Theorem 2.1]{padgett2020analysis}.
\Cref{th:main_i2} of \cref{th:main} follows directly from 
\cref{lem:kern_bd} and \cref{lem:rep1}.

The main result of this article is \cref{cor:final} in \cref{sec:5} below. 
This result provides a complete description of the series representation of the discrete fractional Laplace operator of order $s\in(0,\infty)$. 
The most surprising implication of \cref{cor:final} is that the formula for the function $K \colon \Z \to \R$ 
in \cref{th:main} depends only on the parameter $s \in (0,\infty)$ (cf.\ \cref{def:kern} below). 
In fact, this function is continuous with respect to the parameter $s$ for all $s \in (0,\infty) \backslash\N$ with the points $s\in \N$ all being removable singularities of the function $K \colon \Z \to \R$. 
Hence, we may extend the definition of $K \colon \Z \to \R$ to that of an analytic function (cf.\ \cref{eq:frac_rep_final} of \cref{cor:final}).

The remainder of this article is organized as follows.
In \cref{sec:motivate} we briefly motivate our interest in the development of a series representation 
for the discrete fractional Laplace operator of arbitrary order. 
In particular, we focus on its application to the study of the Anderson localization problem in materials science 
and its application to transport problems in plasma physics. Next, in
\cref{sec:2} we recall several basic definitions and properties of sequence spaces and introduce the 
so-called logarithmic norm. 
Afterwards, in \cref{sec:3} we define the discrete Laplace operator of arbitrary real-valued positive order. 
We do so by introducing the heat semigroup generated by the discrete Laplace operator and then defining 
higher-order powers via induction.
In \cref{sec:4} we define a discrete fractional kernel function and provide a detailed investigation 
of its various quantitative and qualitative properties.
Thereafter, in \cref{sec:5} we construct a series representation for real-valued positive powers of 
the discrete fractional Laplace operator by employing the results from \cref{sec:3,sec:4}.
Finally, in \cref{sec:6}, a number of useful concluding remarks are provided. Continuing avenues of research based on 
the results developed in this article are outlined.


\section{Motivation of study}\label{sec:motivate}

In 1958, P.\ W.\ Anderson suggested that the existence of sufficiently large disorder in a semi-conductor could lead to
spatial localization of electrons \cite{PhysRev.109.1492}.
This localization of electrons in space has since been referred to as \emph{Anderson localization}. 
In an effort to better understand the conditions under which Anderson localization may occur, there have been numerous theoretical and experimental studies of the phenomenon. 
The occurrence of Anderson localization can be defined mathematically through so-called \emph{dynamical localization}, \emph{statistical localization}, or \emph{spectral localization}.
Anderson localization in the dynamical sense is characterized by an exponential decay with respect to time of the wave function 
which represents the particle of interest.
Anderson localization in the statistical sense occurs if the eigenvalues of of the system's associated Hamiltonian are discrete and infinitely close to one 
another when projected onto a finite-dimensional subspace. 
These two definitions of localization have often been related to conventional
strategies used in localized investigations, such as the scaling and perturbation theory (cf., e.g., 
\cite{brandes2003anderson,kramer1993localization,
lagendijk2009fifty,aizenman1993localization,pastur1980spectral,
frohlich1983absence}
and the references therein).
The third definition---the spectral definition---uses the spectrum of the system's infinite-dimensional Hamiltonian operator to 
study the occurrence of Anderson localization (cf., e.g., \cite{MR2215610,MR1779620,MR2219273}).
In particular, this characterization states that, if the Hamiltonian driving the physical system exhibits absolutely continuous spectrum, 
then the system dynamics will not be localized.
The spectral definition becomes the one of interest in our ensuing motivation.

A brief motivation for our interest in the development of a series representation of the discrete fractional Laplace 
operator of arbitrary order is the following.
Let $\mathfrak{C} , s, c, T \in (0,\infty)$, 
let $(\Omega,\mathcal{F},\mathbb{P})$ be a probability space,
let $\ellTwo$ be the set of all $v \colon \Z \to \R$ which satisfy that $\sum_{k\in\Z} \abs{v(k)}^2 < \infty$,
let $\vt{\cdot,\cdot} \colon \ellTwo\times\ellTwo \to \R$ satisfy for all $v,w \in \ellTwo$ that $\vt{v,w} = \sum_{n\in\Z} v(n)w(n)$,
let $K \colon \Z \to \R$ satisfy for all $n\in\Z$ that $K(-n) = K(n)$ and $\abs{K(n)} \le \mathfrak{C} \abs{n}^{-(1+2s)}$, 
let $\varepsilon_n \colon \Omega \to [-\nicefrac{c}{2},\nicefrac{c}{2}]$, $n\in\Z$, be i.i.d.\ random 
variables\footnote{Note that the expression \emph{i.i.d}\ is an abbreviation for the expression 
\emph{independently and identically distributed}}, let $\delta_n \in \ellTwo$, $n\in\Z$, satisfy for all $k,n\in\Z$ 
with $k\neq n$ that $\delta_n(n) = 1$ and $\delta_n(k) = 0$,
and let $u \colon [0,T] \times \Z \times \Omega \to \R$ satisfy for all $t\in[0,T]$, $n\in\Z$ that $u(0,n) = \delta_0(n)$ and
\begin{equation}\label{eq:motivate1}
\pr[\big]{ \tfrac{d}{dt} u}\lrSpace(t,n) = \sum_{k\in\Z} K(k) \pr[\big]{ u(t,n) - u(t,n-k) } + \sum_{k\in\Z} \varepsilon_k \vt{ u(t,k) , \delta_k } \delta_k .
\end{equation}
The situation above is the mathematical formulation of the physical scenario in which electrons are moving through a 
disordered lattice via (possibly) long-range interactions.
The positive value of $c\in(0,\infty)$ represents the maximum magnitude of disorder which can occur at a given point in the lattice $\Z$. 
The i.i.d.\ random variables $\varepsilon_n \colon \Omega \to [-\nicefrac{c}{2},\nicefrac{c}{2}]$, $n\in\Z$, represent the actual disorder at each point in the lattice $\Z$.
Note that the probabilistic nature of this formulation allows for the existence of, e.g., measurement errors.
The function $K \colon \Z \to \R$ describes which long-range jumps are observed to occur as electrons traverse through the lattice $\Z$.
Observe that this implies that the real number $s\in(0,\infty)$ imposes a decay condition on the probability of long-range jumps.
Finally, for every $t\in[0,T]$, $n\in\Z$ it holds that $u(t,n)$ represents the (possibly scaled) probability that the electron will be located at lattice point $n \in \Z$ at time $t \in [0,T]$ (cf., e.g., Padgett et al.\ \cite[Section 2.2]{padgett2020anomalous}).

Based on the above discussion, there arise two immediate questions of interest.
\begin{enumerate}[label=(\Roman*)]
\item\label{q1} What are appropriate (or physically-relevant) choices of the kernel function $K$?
\item\label{q2} Will Anderson localization occur for a system with the long-range interactions described in \cref{eq:motivate1} for all choices of $s,c\in(0,\infty)$?
\end{enumerate}
First, note that there is no unique answer to the question posed in \cref{q1} due to the fact that the construction of mathematical models often depends upon the employed assumptions and individual goals of the scientist constructing them.
It has been recently demonstrated that the discrete fractional Laplace operator is well suited to describe long-range interactions observed in various physical systems, including semi-crystalline polymers and dusty plasma
(cf., e.g., \cite{PhysRevResearch.2.043375,padgett2020anomalous,kostadinova2021active} and the references therein).
In fact, it is the case that when $s \in (0,1)$ one observes so-called \emph{superdiffusion} and when $s\in (1,\infty)$ one observes so-called \emph{subdiffusion}.
Thus, as a starting point we consider the case where the function $K$ coincides with the definition of the discrete fractional Laplace operator.

Next, under the assumption that the function $K$ coincides with the definition of the discrete fractional Laplace operator, we consider the question posed in \cref{q2}.
It is well-known that if $s=1$, then the solution $u$ of \cref{eq:motivate1} will exhibit Anderson localization for all $c \in (0,\infty)$.
What is not known---and an interest which motivates the current study---is whether or not for all $s\in(0,\infty)$ it holds that the solution $u$ of \cref{eq:motivate1} will exhibit Anderson localization for all $c\in(0,\infty)$.
This question can be studied via the following result which follows immediately from Liaw \cite[Corollary 3.2]{liaw2013approach}.

\begin{corollary}\label{cor:motivate}
Let $\mathfrak{C} , s, c, T \in (0,\infty)$, 
let $(\Omega,\mathcal{F},\mathbb{P})$ be a probability space,
let $\ellTwo$ be the set of all $v \colon \Z \to \R$ which satisfy that $\sum_{k\in\Z} \abs{v(k)}^2 < \infty$,
let $\vt{\cdot,\cdot} \colon \ellTwo\times\ellTwo \to \R$ satisfy for all $v,w \in \ellTwo$ that $\vt{v,w} = \sum_{n\in\Z} v(n)w(n)$,
let $K \colon \Z \to \R$ satisfy for all $n\in\Z$ that $K(-n) = K(n)$ and $\abs{K(n)} \le \mathfrak{C} \abs{n}^{-(1+2s)}$, 
let $\varepsilon_n \colon \Omega \to [-\nicefrac{c}{2},\nicefrac{c}{2}]$, $n\in\Z$, be i.i.d.\ random variables, let $\delta_n \in \ellTwo$, $n\in\Z$, satisfy for all $k,n\in\Z$ with $k\neq n$ that $\delta_n(n) = 1$ and $\delta_n(k) = 0$, and
let $H \colon \ellTwo \times \Omega \to \ellTwo$ satisfy for all $u \colon \Z \times \Omega \to \R$, $n\in\Z$ with $\mathbb{P}(u \in \ellTwo) = 1$ that
\begin{equation}\label{eq:motivate2}
(Hu)(n) = \sum_{k\in\Z} K(k) \pr[\big]{ u(n) - u(n-k) } + \sum_{k\in\Z} \varepsilon_k \vt{ u(k) , \delta_k } \delta_k .
\end{equation} 
Then if $H$ has purely singular spectrum\footnote{See, e.g., Kreyszig \cite[Page 371]{kreyszig1978introductory}} $\mathbb{P}$-a.s.\ it holds for all $v\in \ellTwo$ with $\vt{v,v} = 1$ that
\begin{equation}
\mathbb{P}\pr[\bigg]{ \lim_{m\to\infty} \operatorname{dist}\pr[\Big]{ v, \operatorname{span}\cu[\big]{ H^k \delta_0 \colon k \in \{0,1,\ldots,m\} } } = 0 } = 1 .
\end{equation}
\end{corollary}

Observe, that \cref{cor:motivate} shows that in order to study the Anderson localization problem for the operator in \cref{eq:motivate2} above via the spectral definition, we must be able to compute the forward orbit of the operator $H$ with respect to the vector $\delta_0$.
From a numerical perspective, \cref{cor:motivate} implies that if one can find some $v \in \ellTwo$ with $\vt{v,v} = 1$ such that
\begin{equation}
\mathbb{P}\pr[\bigg]{ \lim_{m\to\infty} \operatorname{dist}\pr[\Big]{ v, \operatorname{span}\cu[\big]{ H^k \delta_0 \colon k \in \{0,1,\ldots,m\} } } > 0 } \in (0,1],
\end{equation}
then it holds that $H$ exhibits absolutely continuous spectrum\footnote{See, e.g., Kreyszig \cite[Page 371]{kreyszig1978introductory}} $\mathbb{P}$-a.s.; i.e., Anderson localization does not occur.
Thus, in the case of predicting localization behavior for a system numerically, one would need to be able to compute the action of the operator $H$ \emph{exactly} (or else additional approximations must be introduced).
This need is precisely what motivates our current interest in the development of a series representation of the discrete fractional Laplace operator of arbitrary order.

For improved clarity, we close this section with a few important points. 
The formulation above demonstrates the need to construct an \emph{exact} representation of the \emph{action} of the discrete fractional Laplace operator. 
We accomplish this goal, herein, via the construction of a series representation of the operator (cf.\ \cref{cor:final} below).
The goal of such constructions is predicated on the assumption that the discrete fractional Laplace operator is a good choice for models corresponding to \cref{eq:motivate1}. In fact, the
so-called anomalous diffusion phenomenon has been experimentally observed in various strongly coupled fluids such as ultracold neutral plasma (cf., e.g., Strickler et al.\ \cite{PhysRevX.6.021021}), two-dimensional and quasi-two-dimensional Yukawa liquids (cf., e.g., \cite{PhysRevE.78.026409,PhysRevLett.102.085002,
PhysRevLett.100.055003,PhysRevE.75.016405}), and dusty plasmas (cf., e.g., \cite{PhysRevLett.96.015003,doi:10.1063/1.1449888,
https://doi.org/10.1002/ctpp.200910089}).
As both \emph{subdiffusion} and \emph{superdiffusion} have been observed, it is the case that the discrete 
fractional Laplace operator of arbitrary order is an ideal candidate to model such physical systems, as $s\in(0,1)$ 
may be used to model superdiffusion and $s\in(1,\infty)$ may be used to describe subdiffusion.
The existence of better models for interesting physical systems has been an open question in both mathematics and physics.
We shall leave the study of this issue to our forthcoming papers. 


\section{Background}\label{sec:2}

In this section we review several basic concepts regarding sequence spaces and the logarithmic norm.
More specifically, in \cref{sec:2_1} we introduce the standard $\ellTwo$ sequence space and an 
associated function which we denote the semi-inner product.
In particular, \cref{lem:inner} demonstrates that the standard $\ellTwo$ inner product coincides with our particular semi-inner product.
Afterwards, in \cref{sec:2_2} we define the logarithmic norm and the so-called upper-right Dini derivative.
We then demonstrate a very useful property in \cref{lem:semi_inner} regarding the upper-right Dini derivative of $\ellTwo$ norms.

It is worth noting that the contents of this section have been studied in various parts of the scientific literature 
(although rarely together and in this particular setting).
The concept of semi-inner products has been studied extensively in the literature; cf., e.g., 
\cite{dragomir2004semi,giles1967classes,lumer1961semi}.
They were originally introduced in an effort to extend standard Hilbert space-type arguments to the more 
general setting of normed vector spaces.
Herein, we employ a slight abuse of notation as \cref{def:semi} does not define a semi-inner product in the sense of Lumer (cf., e.g., \cite{lumer1961semi}).
However, the object defined in \cref{def:semi} does possess many of the desired properties of a 
semi-inner product and \cref{lem:inner} demonstrates that no generality is lost by employing this definition.
It is also worth noting that we are not the only authors to employ such notation; cf. e.g.,
S{\"o}derlind
\cite[Definition 5.1]{soderlind2006logarithmic}.
Moreover, it is worth mentioning that \cref{lem:semi_inner} appears in Jones et al.\ \cite[Lemma 2.4]{jones1} in a 
more general setting but we include it below for clarity and completeness.

Throughout this article, $\R$ and $\C$ stand for the usual real and complex number fields, respectively. Further,
let $i = \sqrt{-1}\in\C$, let $\Z$ denote the set of integers, let $\N = \{1,2,3,\ldots\}$, let $\N_0 = \N \cup \{0\}$, 
and for every $z\in\C$ let $\real(z) \in \R$ denote the real part of the complex number $z$.
In addition, we briefly mention a particular notation used throughout this article which emphasizes how various outside 
results are applied. If, for example, a result is referenced which names a particular mathematical object $\mathcal{X}$, 
then in order to state results about a family of objects, herein, e.g., $\mathcal{Y}_t$, $t\in\R$, we will write ``applied for 
every $t\in\R$ with $\mathcal{X} \with \mathcal{Y}_t$ in the notation of \ldots'' in order to clarify its use. We generalize 
this approach in the natural way in the case where multiple mathematical objects are involved (cf., e.g., the proof of 
\cref{lem:lap_semi_props}).
In addition, when carrying out mathematical induction on a variable, say $n \in \N_0$, we will use the notation 
``$\N_0 \ni (n-1) \induct n \in \N$'' to emphasize and clarify both the inductive set and the inductive variable (cf., e.g., the proof of \cref{lem:ver1} below).

\subsection{Sequence spaces}\label{sec:2_1}

\begin{definition}[Set of all sequences]\label{def:seq}
We denote by $\sSpace$ the set of all functions with domain $\Z$ and range $\R$.
\end{definition}

\cfclear
\begin{definition}[The $\ell^2(\Z)$ Hilbert space]\label{def:l2}
\cfconsiderloaded{def:l2}
We denote by $\ellTwo$ the set of all $ u \in \sSpace$ satisfying that $\sum_{k\in\Z} \abs{u(k)}^2 < \infty$ 
\cfload. 
Furthermore, we denote by $\normm{\cdot} \colon \sSpace \to [0,\infty]$ the function which satisfies for all $u \in \sSpace$ that $\normm{u}^2 = \sum_{k\in\Z} \abs{u(k)}^2$\cfout.
\end{definition}

\cfclear
\begin{definition}[The $\ell^2(\Z)$ inner product]\label{def:inner}
\cfconsiderloaded{def:inner}
We denote by $\inner{\cdot,\cdot} \colon \ellTwo \times \ellTwo \to \R$ the function which satisfies for all $u,v\in\ellTwo$ that
$\inner{u,v} = \sum_{k\in\Z} u(k)v(k)$ \cfout.
\end{definition}

\cfclear
\begin{definition}[Semi-inner product]\label{def:semi}
\cfconsiderloaded{def:semi}
We denote by $\semi{\cdot,\cdot} \colon \ellTwo \times \ellTwo \to \R$ the function which satisfies for all $u, v \in \ellTwo$ that
\begin{equation}
\semi{u,v} = \br[\bigg]{ \lim_{\varepsilon\to 0^+} \frac{\normm{v + \varepsilon u} - \normm{v}}{\varepsilon} } \normm{v}
\end{equation}
\cfout.
\end{definition}

\cfclear
\begin{athm}{lemma}{lem:inner}
Let $u,v\in \ellTwo$ \cfload.
Then
\begin{enumerate}[label=(\roman*)]
\item\label{lem:inner_i1} it holds that $\inner{u,u} = \normm{u}^2$ and
\item\label{lem:inner_i2} it holds that
$\inner{u,v} = \semi{u,v}$
\end{enumerate}
\cfout.
\end{athm}

\begin{aproof}
First, observe that \cref{lem:inner_i1} follows immediately from \cref{def:inner,def:l2}. Next, note that \cref{lem:inner_i1} and the fact that $\inner{\cdot,\cdot} \colon \ellTwo \times \ellTwo \to \R$ is a symmetric bilinear form\footnote{It is well known that $\inner{\cdot,\cdot} \colon \ellTwo \times \ellTwo \to \R$ is a symmetric bilinear form as this follows immediately from \cref{def:inner}.}
assure that
\begin{align}
& \semi{u,v} =  \br[\Bigg]{ \lim_{\varepsilon\to 0^+} \frac{\normm{v + \varepsilon u} - \normm{v} }{\varepsilon} } \normm{v}
= \br[\Bigg]{ \lim_{\varepsilon\to 0^+} \frac{\normm{v + \varepsilon u} - \normm{v} }{\varepsilon} \cdot \frac{ \normm{v + \varepsilon u} + \normm{v} }{ \normm{v + \varepsilon u} + \normm{v} } } \normm{v} \nonumber \\
& \quad = \br[\Bigg]{ \lim_{\varepsilon\to 0^+} \frac{\normm{v + \varepsilon u}^2 - \normm{v}^2 }{\varepsilon \pr[\big]{ \normm{v + \varepsilon u} + \normm{v} }} } \normm{v}
= \br[\Bigg]{ \lim_{\varepsilon\to 0^+} \frac{\inner{v + \varepsilon u,v+\varepsilon u} - \inner{v,v} }{\varepsilon \pr[\big]{ \normm{v + \varepsilon u} + \normm{v} }} } \normm{v} \\
& \quad = \br[\Bigg]{ \lim_{\varepsilon\to 0^+} \frac{ \inner{v,v} + 2 \varepsilon \inner{u,v} + \varepsilon^2 \inner{u,u} - \inner{v,v} }{\varepsilon \pr[\big]{ \normm{v + \varepsilon u} + \normm{v} }} } \normm{v}
= \br[\Bigg]{ \lim_{\varepsilon\to 0^+} \frac{2 \varepsilon \inner{u,v} + \varepsilon^2 \inner{u,u} }{\varepsilon \pr[\big]{ \normm{v + \varepsilon u} + \normm{v} }} } \normm{v} \nonumber \\
& \quad = \br[\Bigg]{ \lim_{\varepsilon\to 0^+} \frac{2  \inner{u,v} + \varepsilon \inner{u,u} }{\normm{v + \varepsilon u} + \normm{v} } } \normm{v}
= \br[\Bigg]{ \frac{ 2 \inner{u,v} }{ 2 \normm{v} }  } \normm{v} = \inner{u,v} \nonumber
\end{align}
\cfout.
This establishes \cref{lem:inner_i2}.
\end{aproof}

\subsection{The logarithmic norm}\label{sec:2_2}

\cfclear
\begin{definition}[Logarithmic norm]\label{def:log}
\cfconsiderloaded{def:log}
For every $A \colon \ellTwo \to \ellTwo$ we denote by $\logNorm{A} \in \R$ the real number which satisfies that
\begin{equation}
\logNorm{A} = \sup_{ \substack{ v \in \ellTwo \\ \normm{v} \neq 0 } } \frac{ \inner{ Av, v } }{ \normm{v}^2 } \ifnocf.
\end{equation}
\cfout.
\end{definition}

\cfclear
\begin{definition}[Upper-right Dini derivative]\label{def:dini}
\cfconsiderloaded{def:dini}
For every $v\colon [0,\infty) \to \R$ we denote by $\dini{t} v(t) \in [-\infty,\infty]$, $t\in[0,\infty)$, the function which satisfies for all $t\in[0,\infty)$ that
\begin{equation}
\dini{t} v(t) = \limsup_{\varepsilon\to 0^+} \frac{v(t+\varepsilon) - v(t)}{\varepsilon} .
\end{equation}
\end{definition}

\cfclear
\begin{athm}{lemma}{lem:semi_inner}
It holds for all $t\in[0,\infty)$ and differentiable $v \colon [0,\infty) \to \ellTwo$ that
\begin{equation}
\dini{t} \normm{v(t)} = \br[\bigg]{ \frac{ \inner{\frac{d}{dt} v(t), v(t) } }{ \normm{v(t)}^2 } } \normm{v(t)} \ifnocf.
\end{equation}
\cfout.
\end{athm}

\begin{aproof}
Throughout this proof let $v \colon [0,\infty) \to \ellTwo$, let $t \in[0,\infty)$, and assume without loss of generality that $\normm{v(t)} \neq 0$ \cfload.
Note that the hypothesis that $v$ is differentiable and Taylor's theorem (cf., e.g., Cartan et al.\ \cite[Theorem 5.6.3]{cartan2017differential}) yield that there exist $\delta_t(\varepsilon) \in \ellTwo$, $\varepsilon \in \R$, such that for all $\varepsilon \in \R$ with $\abs{\varepsilon}$ sufficiently small it holds that
\begin{enumerate}[label=(\Alph*)]
\item\label{aa} $v(t+\varepsilon) = v(t) + \varepsilon \frac{d}{dt} v(t) + \abs{\varepsilon} \delta_t(\varepsilon)$ and
\item\label{bb} $\lim_{\varepsilon\to 0} \delta_t(\varepsilon) = 0$.
\end{enumerate}
Combining \cref{aa,bb} with \cref{lem:inner_i2} of \cref{lem:inner}
hence shows that
\begin{align}
& \dini{t} \normm{v(t)} = \limsup_{\varepsilon \to 0^+} \frac{\normm{v(t) + \varepsilon \frac{d}{dt} v(t) + \abs{\varepsilon} \delta_t(\varepsilon)} - \normm{v(t)}}{\varepsilon} \nonumber \\
& \quad = \limsup_{\varepsilon \to 0^+} \frac{\normm{v(t) + \varepsilon \frac{d}{dt} v(t)} - \normm{v(t)}}{\varepsilon}
= \lim_{\varepsilon \to 0^+} \frac{\normm{v(t) + \varepsilon \frac{d}{dt} v(t)} - \normm{v(t)}}{\varepsilon}  \\
& \quad = \br[\bigg]{ \lim_{\varepsilon \to 0^+} \frac{\normm{v(t) + \varepsilon \frac{d}{dt} v(t)} - \normm{v(t)}}{\varepsilon} } \frac{\normm{v(t)}^2}{\normm{v(t)}^2}
= \frac{ \inner[\big]{\frac{d}{dt}v(t), v(t)} }{\normm{v(t)}^2} \normm{v(t)} \nonumber
\end{align}
\cfload.
\end{aproof}

\cfclear
We close \cref{sec:2_2} with a brief discussion of \cref{def:log}.
For every $A \colon \ellTwo \to \ellTwo$ let $\norm{A}_{\text{op}} = \inf\{ c \in [0,\infty] \colon \forall\, v \in \ellTwo \text{ it holds that } \normm{Av} \le c \normm{v} \}$ \cfload. Then \cref{def:log,def:semi} and the Rayleigh quotient theorem (cf., e.g., \cite[Theorem A.26]{driverfunctional}) imply that for every $A \colon \ellTwo \to \ellTwo$ with $\norm{A}_{\text{op}} \in [0,\infty)$, $A = A^*$, and nonempty pure point spectrum\footnote{See, e.g., Kreyszig \cite[Page 371]{kreyszig1978introductory}} it holds that
\begin{equation}
\logNorm{A} = \sup_{ \substack{ v \in \ellTwo \\ \normm{v} \neq 0 } } \frac{ \inner{ Av, v } }{ \normm{v}^2 }
= \max\{ \lambda \in \R \colon \exists \, v \in \ellTwo \text{ with } \normm{v} \neq 0 \text{ and } Av = \lambda v \}
\end{equation}
(e.g., $\logNorm{A}$ is the maximal eigenvalue of $A$).
This fact will prove quite useful in the proof of \cref{lem:lap_semi_props} in \cref{sec:3_1}.


\section{The discrete Laplace operator of arbitrary order}\label{sec:3}

\cfclear
In this section we introduce the discrete fractional Laplace operator and define the notion of real-valued positive powers of this operator.
First, in \cref{sec:3_1} we define the discrete Laplace operator as well as introduce and study its associated discrete heat semigroup. 
\cref{prop:high_lap_defined} is presented in order to clarify the fact that positive integer powers of the discrete Laplace operator map elements of $\ellTwo$ into $\ellTwo$ \cfload.
The associated discrete heat semigroup is shown to be a strongly continuous contraction semigroup via the tools developed in \cref{sec:2_2}.
Note that \cref{def:gamma} is provided for clarity, as the evaluation of the Gamma function with arguments whose real parts are negative occurs frequently throughout the remainder of this article.

Afterwards, in \cref{sec:3_2} we provide a series representation of positive integer powers of the discrete Laplace operator.
The result in \cref{lem:ver1} is well-known in the literature, but its proof is included for completeness (cf., e.g., Kelley and Peterson \cite[Eq.\ (2.1)]{kelley2001difference}).
The series representation presented in \cref{lem:ver1} will be a crucial component in proving the main result of this article (cf.\ \cref{cor:final}).

In \cref{sec:3_3} we define arbitrary real-valued positive powers of the discrete Laplace operator (cf.\ \cref{def:frac_lap}). 
This is accomplished by first defining the case when the positive real-valued power is bounded above by one (cf.\ \cref{def:frac_lap_pre}).
We then define higher-order positive real-valued powers via an inductive procedure.
This definition is shown to be well-defined in $\ellTwo$ in  \cref{lem:background}; i.e., it is shown that the discrete fractional Laplace operator maps $\ellTwo$ into $\ellTwo$.
At this point, we wish to again emphasize that \cref{def:frac_lap} is not a direct ``discretization'' of the pointwise formula for the continuous case (cf., e.g., \cref{eq:caff} for the case where $s \in (0,1)$), but rather the fractional power of the discrete Laplace operator.

\subsection{The discrete Laplace operator and its associated semigroup}\label{sec:3_1}

\cfclear
\begin{definition}[Gamma function]\label{def:gamma}
\cfconsiderloaded{def:gamma}
\cfload
Let $X = \{ z \in \C \colon \real(z) \in (0,\infty) \}$ and let $\tilde{\Gamma} \colon X \to \C$ be the function which satisfies for all $z \in X$ that $\tilde{\Gamma}(z) = \int_0^\infty x^{z-1} \exp(-x) \, dx$.
Then we denote\footnote{Note that $\floor{\cdot} \colon \Z \to \R$ satisfies for all $x \in \R$ that $\floor{x} = 
\max\{ n \in \Z \colon n \le x \}$.} 
by $\gamfn \colon \C \to \C$ the function which satisfies for all $ v \in X $, $w,z\in\C$ with $\real(z) \in (-\infty,0]\backslash\{\ldots,-2,-1,0\}$ and $\real(w) \in \{\ldots,-2,-1,0\}$ that
$\gamfn(v) = \tilde{\Gamma}(v)$, $\nicefrac{1}{\gamfn(w)} = 0$,  and
\begin{equation}
\gamfn(z) = 
\frac{\tilde{\Gamma}(z + \abs{\floor{\real(z)}})}{(z+\abs{\floor{\real(z)}}-1) (z+\abs{\floor{\real(z)}}-2) \cdot \ldots \cdot z }.
\end{equation}
\end{definition}

\cfclear
\begin{athm}{proposition}{prop:high_lap_defined}
It holds\footnote{Note that for all $n\in\N$, $k \in \{0,1,\ldots,n\}$ it holds that $\binom{n}{k} = \nicefrac{\gamfn(1+n)}{(\gamfn(1+k) \gamfn(1+n-k))}$ (cf.\ \cref{def:gamma}).}
for all $s \in \N$, $u \in \ellTwo$ that
\begin{equation}
\sum_{n\in\Z} \left\lvert \sum_{k=0}^{2s} (-1)^{k-s} \binom{2s}{k} u(n-s+k) \right\rvert^2 < \infty
\end{equation}
\cfload.
\end{athm}

\begin{aproof}
Throughout this proof let $s\in\N$, $u \in \ellTwo$. Observe that the fact that $u \in \ellTwo$ ensures that $\sum_{k\in\Z} \abs{u(k)}^2 < \infty$. This, the triangle inequality, the fact that $s\in\N$, and the fact Jensen's inequality implies that for all $r,m\in\N$, $v_1,v_2,\ldots,v_m\in[0,\infty)$ it holds that $[\sum_{k=1}^m v_k]^r \le m^{\max\{r-1,0\}} \sum_{k=1}^m v_k^r$ assure that
\begin{align}\label{higher_bd}
& \sum_{n\in\Z} \left\lvert \sum_{k=0}^{2s} (-1)^{k-s} \binom{2s}{k} u(n-s+k) \right\rvert^2
\le \sum_{n\in\Z} \left[ \sum_{k=0}^{2s} \left\lvert \binom{2s}{k} u(n-s+k) \right\rvert \right]^2 \\
& \quad \le \sum_{n\in\Z} \left[ 2s \sum_{k=0}^{2s} \binom{2s}{k} \left\lvert u(n-s+k) \right\rvert^2 \right]
= 2s \sum_{k=0}^{2s} \binom{2s}{k} \left[ \sum_{n\in\Z} \abs{u(n-s+k)}^2 \right] < \infty. \nonumber
\end{align}
\end{aproof}

\cfclear
\begin{definition}[Discrete Laplace operator]\label{def:lap}
\cfconsiderloaded{def:lap}
\cfload
We denote by $\lap \colon \ellTwo \to \ellTwo$ the function which satisfies for all $u \in \ellTwo$, $n\in\Z$ that
\begin{equation}
(\lap u)(n) = u(n-1) - 2u(n) + u(n+1) \ifnocf.
\end{equation}
\cfout[.]
\end{definition}

\cfclear
\begin{definition}[Identity operator]\label{def:id}
\cfconsiderloaded{def:id}
\cfload
We denote by $\id \colon \sSpace \to \sSpace$ the operator which satisfies for all $u\in\sSpace$, $n\in\Z$
that $(\id u)(n) = u(n)$ \cfout.
\end{definition}

\cfclear
\begin{athm}{proposition}{prop:semi_defined}
Let $u \in \ellTwo$ \cfload.
Then it holds for all $z\in[0,\infty)$ that
\begin{equation}
\sum_{n\in\Z} \left \lvert \sum_{k\in\Z} \exp( -2z ) \br[\Bigg]{ \sum_{j\in\N_0} \frac{ z^{2j + n - k} }{ \gamfn( 1+j ) \gamfn( j + n - k + 1 ) } } u(k) \right \rvert^2  < \infty
\end{equation}
\cfout.
\end{athm}

\begin{aproof}
Throughout this proof let $I_k \colon [0,\infty) \to \R$, $k\in\Z$, satisfy for all $z\in[0,\infty)$, $k\in\Z$ that
\begin{equation}\label{eq:bessel}
I_k(z) = \sum_{j\in\N_0} \frac{ \pr[\big]{ \nicefrac{z}{2} }^{2j + k} }{ \gamfn( 1+j ) \gamfn( j + k + 1 ) }
\end{equation} 
\cfload.
Observe that \cref{eq:bessel} and Olver et al.\ \cite[Eq.\ 10.27.1]{olver2010nist}
(applied for every $k \in \Z$ with $I_k \with I_k$, $z \with 2z$ in the notation of Olver et al.\ \cite[Eq.\ 10.27.1]{olver2010nist})
assure that for all $z \in [0,\infty)$ it holds that
\begin{align}
& \sum_{k\in\Z} \br[\Bigg]{ \sum_{j\in\N_0} \frac{ z^{2j + k} }{ \gamfn( 1+j ) \gamfn( j + k + 1 ) } } 
= \sum_{k\in\Z} I_k(2z) \\
& \quad = I_0(2z)  + \sum_{k\in\N} \pr[\big]{ I_k(2z) + I_{-k}(2z) } 
= I_0(2z) + 2 \sum_{k\in\N} I_k(2z)
. \nonumber
\end{align}
Combining this, \cref{eq:bessel}, 
the triangle inequality,
Minkowski's inequality,
and
Olver et al.\ \cite[Eq.\ 10.35.5]{olver2010nist} (applied with $I_k \with I_k$, $z \with 2z$ in the notation of Olver et al.\ \cite[Eq.\ 10.35.5]{olver2010nist}) ensures that for all $z\in [0,\infty)$ it holds that
\begin{align}
& \sum_{n\in\Z} \left \lvert \sum_{k\in\Z} \exp( -2z ) \br[\Bigg]{ \sum_{j\in\N_0} \frac{ z^{2j + n - k} }{ \gamfn( 1+j ) \gamfn( j + n - k + 1 ) } } u(k) \right \rvert^2 \nonumber \\
& \quad = \sum_{n\in\Z} \left \lvert \sum_{k\in\Z} \exp(-2z) I_{k}(2z) u(n-k) \right \rvert^2 
\le \exp(-4z) \left[ \sum_{k\in\Z} \abs{I_{k}(2z)} \left( \sum_{n\in\Z} \abs{ u(n-k) }^2 \right)^{\!\!\nicefrac{1}{2}} \right]^{\!2} \nonumber \\
& \quad = \exp(-4z)\left[ \sum_{k\in\Z} I_{k}(2z) \normm{u} \right]^{\!2} 
= \exp(-4z) \left[ I_0(2z) + 2 \sum_{k\in\N} I_k(2z) \right]^{\!2} \normm{u}^2 \\
& \quad = \exp(-4z) \br[\big]{ \exp(2z) }^2 \normm{u}^2
= \normm{u}^2 < \infty. \nonumber
\end{align}
\end{aproof}

\cfclear
\begin{definition}[Discrete heat semigroup]\label{def:semigroup}
\cfload
We denote by $\semigroup_z(\lap) \colon \ellTwo \to \ellTwo$, $z\in[0,\infty)$, the function which satisfies for all $z\in[0,\infty)$, $u \in \ellTwo$, $n\in\Z$ that
\begin{equation}\label{this_semi}
\pr[\big]{ \semigroup_z(\lap) u }(n) = \sum_{k\in\Z} \exp( -2z ) \br[\Bigg]{ \sum_{j\in\N_0} \frac{ z^{2j + n - k} }{ \gamfn( 1+j ) \gamfn( j + n-k + 1 ) } } u(k)\ifnocf.
\end{equation}
\cfout[.]
\end{definition}

\cfclear
\begin{athm}{lemma}{lem:lap_semi_props}
Let $ u \in \ellTwo$ \cfload.
Then
\begin{enumerate}[label=(\roman*)]
\item\label{lem:lap_semi_props_i1} it holds that $\semigroup_z(\lap) \colon \ellTwo \to \ellTwo$, $z\in[0,\infty)$, is a strongly continuous semigroup\footnote{Cf., e.g., Jones et al.\ \cite[Definition 2.6]{jones1}},
\item\label{lem:lap_semi_props_i1a} it holds for all $z\in[0,\infty)$, $n\in\Z$ that $\frac{d}{dz} (\semigroup_z(\lap) u)(n) = (\lap \semigroup_z(\lap) u)(n)$,
\item\label{lem:lap_semi_props_i2} it holds that $\logNorm{\lap} \in (-\infty,0)$, 
and
\item\label{lem:lap_semi_props_i3}
it holds for all $z \in [0,\infty)$ that $\normm{\semigroup_z(\lap) u } \le \exp(z \logNorm{\lap}) \normm{u} \le \normm{u}$
\end{enumerate}
\cfout.
\end{athm}

\begin{aproof}
First, note that $\semigroup_z(\lap) \colon \ellTwo \to \ellTwo$, $z\in[0,\infty)$, is a strongly continuous semigroup if it holds for all $t,z \in [0,\infty)$, $v \in \ellTwo$ that
\begin{enumerate}[label=(\Alph*)]
\item\label{aaa} $\semigroup_0(\lap) = \id$,
\item\label{bbb} $\semigroup_{t+z}(\lap) = \semigroup_t(\lap) \semigroup_z(\lap)$, and
\item\label{ccc} $\lim_{z\to 0^+} \normm{ \semigroup_z(\lap) v - v } = 0$
\end{enumerate}
\cfload.
Observe that the fact that for all $v \in \ellTwo$ it holds that $\sup_{n\in\Z} \abs{v(n)} < \infty$
and, e.g., Ciaurri et al.\ \cite[Proposition 1]{ciaurri2017harmonic} 
(applied for every $v \in \ellTwo$ with $ f \with v $, $W_t \with \semigroup_z(\lap)$ in the notation of Ciaurri et al.\ \cite[Proposition 1]{ciaurri2017harmonic})
show that $\semigroup_z(\lap) \colon \ellTwo \to \ellTwo$, $z\in[0,\infty)$ satisfies \cref{aaa,bbb,ccc}.
This establishes \cref{lem:lap_semi_props_i1}.
Next, note that, e.g., Ciaurri et al.\ \cite[Proposition 2]{ciaurri2017harmonic} 
(applied for every $n \in \Z$ with $ f \with u $, $(u(n,t))_{t\in[0,\infty)} \with ((\semigroup_z(\lap)u)(n))_{z\in[0,\infty)}$ in the notation of Ciaurri et al.\ \cite[Proposition 2]{ciaurri2017harmonic}) establishes \cref{lem:lap_semi_props_i1a}. 
In addition, observe that if for all $v \in \ellTwo$ with $\normm{v} \neq 0$ it holds that $\inner{\lap v, v} \in (-\infty,0)$ then it holds that $\logNorm{\lap} \in (-\infty,0)$ \cfload.
To that end, note that 
for all $v \in \ellTwo$ it holds that
\begin{align}
\inner{ \lap v , v } & = \sum_{k\in\Z} \br[\big]{ (\lap v)(k) }  v(k) 
= \sum_{k\in\Z} \pr[\big]{ v(k-1) - 2v(k) + v(k+1) } v(k) \nonumber \\
& = \sum_{k\in\Z} v(k-1)v(k) - 2\sum_{k\in\Z} v(k)^2 + \sum_{k\in\Z} v(k+1)v(k) \nonumber \\
& = \sum_{k\in\Z} v(k-1)v(k) - \sum_{k\in\Z} v(k)^2 - \sum_{k\in\Z} v(k-1)^2 + \sum_{k\in\Z} v(k)v(k-1) \\
& = - \sum_{k\in\Z} \br[\big]{ v(k-1)^2 - 2v(k-1)v(k) + v(k)^2 }
= - \sum_{k\in\Z} \br[\big]{ v(k-1) - v(k) }^2 . \nonumber
\end{align}
This demonstrates that for all $v\in\ellTwo$ with $\normm{v} \neq 0$ it holds that $\inner{\lap v, v} \in (-\infty,0)$. This establishes \cref{lem:lap_semi_props_i2}.
Moreover, observe that \cref{lem:lap_semi_props_i1}, \cref{lem:lap_semi_props_i2}, and Jones et al. \cite[Lemma 2.8]{jones1} (applied with $x \with u$, $A \with \lap$, $(\mathbb{T}_t(A))_{t\in[0,\infty)} \with (\semigroup_z(\lap))_{z\in[0,\infty)}$ in the notation of Jones et al. \cite[Lemma 2.8]{jones1}) hence prove \cref{lem:lap_semi_props_i3}.
\end{aproof}

\subsection{Positive integer powers of the discrete Laplace operator}\label{sec:3_2}

\cfclear
\begin{athm}{lemma}{lem:ver1}
\cfconsiderloaded{def:id}
\cfload
It holds
for all $s \in \N$, $u \in \ellTwo$, $n \in \Z$ that
\begin{equation}\label{eq:ver1_1}
\pr[\big]{ \fracLap u }(n) = \sum_{k=0}^{2s} (-1)^{k-s} \binom{2s}{k} u(n-s+k) \ifnocf.
\end{equation}
\cfload.
\end{athm}

\begin{aproof}
We prove \cref{eq:ver1_1} by induction on $s\in\N$. For the base case $s=1$ observe that \cref{def:lap} establishes \cref{eq:ver1_1}. This proves \cref{eq:ver1_1} in the case $s=1$. For the induction step $\N \ni (s-1) \induct s \in \N \cap [2,\infty)$, let $s \in \N \cap [2,\infty)$ and assume for all $\mathfrak{s} \in \{1,2,\dots,s-1\}$, $u \in \ellTwo$, $n\in\Z$ that \cref{eq:ver1_1} holds. Note that 
the fact that $s \in \N \cap [2,\infty)$ implies that for all $u \in \ellTwo$, $n\in\Z$ it holds that 
\begin{equation}\label{eq:semi_var}
\pr[\big]{ \fracLap u }(n) = \pr[\big]{ -\lap ( (-\Delta)^{s-1} u ) }(n).
\end{equation}
This and the induction hypothesis ensure that for all $u \in \ellTwo$, $n\in\Z$ it holds that
\begin{align}
\pr[\big]{ \fracLap u }(n) & = \pr[\bigg]{ -\lap \pr[\bigg]{ \SmallSum_{k=0}^{2(s-1)} (-1)^{k-(s-1)} \binom{2(s-1)}{k} u(\cdot-(s-1)+k) } }(n) \nonumber \\
& = 2 \SmallSum_{k=0}^{2s-2} (-1)^{k-s+1} \Binom{2s-2}{k} u(n-s+1+k) \nonumber \\
& \qquad - \SmallSum_{k=0}^{2s-2} (-1)^{k-s+1} \Binom{2s-2}{k} u((n-1)-s+1+k) \nonumber \\
& \qquad - \SmallSum_{k=0}^{2s-2} (-1)^{k-s+1} \binom{2s-2}{k} u((n+1)-s+1+k) \nonumber \\
& = 2 \SmallSum_{k=0}^{2s-2} (-1)^{k-s+1} \Binom{2s-2}{k} u(n-s+1+k) - \SmallSum_{k=0}^{2s-2} (-1)^{k-s+1} \Binom{2s-2}{k} u(n-s+k) \nonumber \\
& \qquad - \SmallSum_{k=0}^{2s-2} (-1)^{k-s+1} \Binom{2s-2}{k} u(n+2-s+k) \\
& = 2 \SmallSum_{k=1}^{2s-1} (-1)^{k-s} \Binom{2s-2}{k-1} u(n-s+k) + \SmallSum_{k=0}^{2s-2} (-1)^{k-s} \Binom{2s-2}{k} u(n-s+k) \nonumber \\
& \qquad + \SmallSum_{k=2}^{2s} (-1)^{k-s} \Binom{2s-2}{k-2} u(n-s+k) \nonumber \\
& = \br[\bigg]{ \SmallSum_{k=1}^{2s-1} (-1)^{k-s} \Binom{2s-2}{k-1} u(n-s+k) + \SmallSum_{k=0}^{2s-2} (-1)^{k-s} \Binom{2s-2}{k} u(n-s+k) } \nonumber \\
& \qquad + \br[\bigg]{ \SmallSum_{k=1}^{2s-1} (-1)^{k-s} \Binom{2s-2}{k-1} u(n-s+k) + \SmallSum_{k=2}^{2s} (-1)^{k-s} \Binom{2s-2}{k-2} u(n-s+k) }. \nonumber
\end{align}
Combining this and the fact that for all $n\in\N$, $k \in \{1,2,\dots,n-1\}$ it holds that $\binom{n}{k} = \binom{n-1}{k} + \binom{n-1}{k-1}$ assures that for all $u \in \ellTwo$, $n\in\Z$ it holds that
\begin{align}\label{eq:b3}
\pr[\big]{ \fracLap u }(n)
& = (-1)^{s-1} u(n+s-1) + (-1)^{-s} u(n-s) + \SmallSum_{k=1}^{2s-2} (-1)^{k-s} \Binom{2s-1}{k} u(n-s+k) \nonumber \\
& \qquad + (-1)^{1-s} u(n+1-s) + (-1)^{s} u(n+s) \nonumber \\
& \qquad + \SmallSum_{k=2}^{2s-1} (-1)^{k-s} \Binom{2s-1}{k-1} u(n-s+k) \nonumber \\
& = (-1)^{s-1} u(n+s-1) + (-1)^{-s} u(n-s) + (-1)^{1-s} u(n+1-s) \nonumber \\
& \qquad + (-1)^{s} u(n+s) + (-1)^{1-s} \textstyle\Binom{2s-1}{1} u(n-s+1) \nonumber \\
& \qquad + (-1)^{s-1} \textstyle\Binom{2s-1}{2s-2} u(n-s+(2s-1)) \\
& \qquad + \SmallSum_{k=2}^{2s-2} (-1)^{k-s} \br[\Big]{ \Binom{2s-1}{k} + \Binom{2s-1}{k-1} } u(n-s+k) \nonumber \\
& = (-1)^{s-1} u(n+s-1) + (-1)^{-s} u(n-s) + (-1)^{1-s} u(n+1-s) \nonumber \\
& \qquad + (-1)^{s} u(n+s) + (-1)^{1-s} \Binom{2s-1}{1} u(n-s+1) \nonumber \\
& \qquad + (-1)^{s-1} \Binom{2s-1}{2s-2} u(n+s-1) + \SmallSum_{k=2}^{2s-2} (-1)^{k-s} \Binom{2s}{k} u(n-s+k). \nonumber
\end{align}
Next, observe that the fact that for all $n\in\N$, $k \in \{1,2,\dots,n-1\}$ it holds that $\binom{n}{k} = \binom{n-1}{k} + \binom{n-1}{k-1}$ and the fact that for all $n\in\N$, $k \in \{0,1,\dots,n\}$ it holds that $\binom{n}{k} = \binom{n}{n-k}$ show that
\begin{equation}\label{eq:b1}
1 + \Binom{2s-1}{2s-2} = \Binom{2s-1}{2s-1} + \Binom{2s-1}{2s-2} = \Binom{2s}{2s-1}
\end{equation}
and
\begin{equation}\label{eq:b2}
1 + \Binom{2s-1}{1} = \Binom{2s-1}{0} + \Binom{2s-1}{1} = \Binom{2s}{1}. 
\end{equation}
Combining \cref{eq:b1,eq:b2,eq:b3} hence implies that for all $u \in \ellTwo$, $n\in\Z$ it holds that
\begin{align}
\pr[\big]{ \fracLap u }(n)
& = (-1)^{-s} \Binom{2s}{0} u(n-s) + (-1)^{s} \Binom{2s}{2s} u(n+s) + (-1)^{1-s} \Binom{2s}{1} u(n-s+1) \nonumber \\
& \qquad + (-1)^{s-1} \Binom{2s}{2s-1} u(n+s-1) + \SmallSum_{k=2}^{2s-2} (-1)^{k-s} \Binom{2s}{k} u(n-s+k) \\
& = \SmallSum_{k=0}^{2s} (-1)^{k-s} \Binom{2s}{k} u(n-s+k). \nonumber
\end{align}
Induction therefore establishes \cref{eq:ver1_1}. 
\end{aproof}


\subsection{The discrete fractional Laplace operator of arbitrary order}\label{sec:3_3}

\cfclear
\begin{athm}{lemma}{lem:background}
Let $m \in \N$, $s\in(m-1,m)$.
Then
it holds for all $u \in \ellTwo$ that
\begin{equation}\label{eq:2.14_int}
\normm[\bigg]{ \frac{1}{\gamfn(-(s-m+1))} \int_0^\infty z^{-s+m-2} \br[\big]{ \semigroup_z(\lap) - \id } \pr[\big]{ (-\lap)^{m-1} u } \, dz } < \infty
\end{equation}
\cfout.
\end{athm}

\begin{aproof}
Throughout this proof let $ m \in \N$, $s\in (m-1,m)$.
We claim that for all $z\in (0,\infty)$, $u \in \ellTwo$ it holds that
\begin{equation}\label{eq:2.15_claim}
\normm[\big]{ ( \semigroup_z(\lap) - \id ) u } \le \int_0^z \exp\pr[\big]{ (z-w) \logNorm{\lap} } \normm{\lap u} \, dw \ifnocf.
\end{equation}
\cfload.
Note that the fact that $\inner{\cdot,\cdot} \colon \ellTwo \times \ellTwo \to \R$ is a symmetric bilinear form ensures that for all $ u \in \ellTwo $, $z\in(0,\infty)$ it holds that
\begin{align}
& \frac{d}{dz} \normm[\big]{ \pr[]{ \semigroup_z(\lap) - \id } u } = \frac{ \inner{ \frac{d}{dz} (\semigroup_z(\lap) - \id) u , (\semigroup_z(\lap) - \id )u } }{ \normm{ (\semigroup_z(\lap) - \id ) u } } \\
& \quad = \frac{ \inner{\lap \semigroup_z(\lap) u , (\semigroup_z(\lap) - \id )u } }{ \normm{ (\semigroup_z(\lap) - \id ) u } } 
= \frac{ \inner{\lap ( \semigroup_z(\lap) - \id ) u , (\semigroup_z(\lap) - \id )u } }{ \normm{ (\semigroup_z(\lap) - \id ) u } } + \frac{ \inner{\lap u , (\semigroup_z(\lap) - \id )u } }{ \normm{ (\semigroup_z(\lap) - \id ) u } } \nonumber \\
& \quad \le \br[\Bigg]{ \sup_{v \in \ellTwo} \frac{ \inner{\lap v , v } }{ \normm{ v } } } + \frac{ \inner{\lap u , (\semigroup_z(\lap) - \id )u } }{ \normm{ (\semigroup_z(\lap) - \id ) u } }
= \logNorm{\lap} + \frac{ \inner{\lap u , (\semigroup_z(\lap) - \id )u } }{ \normm{ (\semigroup_z(\lap) - \id ) u } } \ifnocf. \nonumber
\end{align}
\cfload.
This, the Cauchy-Swartz inequality, and the fact that \cref{lem:lap_semi_props_i1} of \cref{lem:lap_semi_props} implies that for all $u \in \ellTwo$ it holds that 
$\lim_{z\to 0^+} \normm{ (\semigroup_z(\lap) \allowbreak - \allowbreak \id )u } = 0$
assure that for all $z\in(0,\infty)$, $u \in \ellTwo$ it holds that
\begin{align}
& \normm[\big]{ \pr[]{ \semigroup_z(\lap) - \id } u }
\le \int_0^z \exp\pr[\big]{ (z-w) \logNorm{\lap} } \frac{ \abs[\big]{ \inner{\lap u , (\semigroup_w(\lap) - \id )u }  }}{ \norm{ (\semigroup_w(\lap) - \id ) u }_2 } \, dw \\
& \le \int_0^z \exp\pr[\big]{ (z-w) \logNorm{\lap} } \frac{ \normm{\lap u} \normm{ (\semigroup_w(\lap) - \id )u } }{ \normm{ (\semigroup_w(\lap) - \id ) u } } \, dw
= \int_0^z \exp\pr[\big]{ (z-w) \logNorm{\lap} } \normm{\lap u} \, dw. \nonumber
\end{align}
Combining this and the fact that \cref{prop:high_lap_defined} (applied with $s \with 1$ in the notation of \cref{prop:high_lap_defined}) ensures that for all $u \in \ellTwo$ it holds that $\normm{\lap u} < \infty$ proves \cref{eq:2.15_claim}.
Next, observe that \cref{eq:2.15_claim}, \cref{lem:ver1}, and Jensen's inequality guarantee that for all $u \in \ellTwo$ it holds that
\begin{align}\label{eq:2.18_int}
& \normm[\bigg]{ \frac{1}{\gamfn(-(s-m+1))} \int_0^\infty z^{-s+m-2} \br[\big]{ \semigroup_z(\lap) - \id } \pr[\big]{ (-\lap)^{m-1} u } \, dz } \nonumber \\
& \quad \le \frac{1}{\abs{\gamfn(-(s-m+1))}} \int_0^\infty z^{-s+m-2} \normm[\big]{ \br[\big]{ \semigroup_z(\lap) - \id } \pr[\big]{ (-\lap)^{m-1} u } } \, dz \\
& \quad \le \frac{\normm{(-\lap)^m u}}{\abs{\gamfn(-(s-m+1))}} \int_0^\infty z^{-s+m-2} \br[\bigg]{ \int_0^z \exp\pr[\big]{ (z-w) \logNorm{\lap} } \, dw } \, dz \ifnocf. \nonumber 
\end{align}
\cfload.
In addition, note that \cref{lem:lap_semi_props_i2} of \cref{lem:lap_semi_props}
and integration by parts show that
\begin{align}\label{eq:2.19_int}
0 & \le \int_0^\infty z^{-s+m-2} \br[\bigg]{ \frac{ 1 - \exp\pr[\big]{ z \logNorm{\lap} } }{ -\logNorm{\lap} } } \, dz \nonumber \\
& = \int_0^1 z^{-s+m-2} \br[\bigg]{ \frac{ 1 - \exp\pr[\big]{ z \logNorm{\lap} } }{ -\logNorm{\lap} } } \, dz + \int_1^\infty z^{-s+m-2} \br[\bigg]{ \frac{ 1 - \exp\pr[\big]{ z \logNorm{\lap} } }{ -\logNorm{\lap} } } \, dz \nonumber \\
& \le \int_0^1 z^{-s+m-2} \br[\bigg]{ \frac{ 1 - \exp\pr[\big]{ z \logNorm{\lap} } }{ -\logNorm{\lap} } } \, dz + \int_1^\infty z^{-s+m-2} \, dz \nonumber \\
& = \lim_{w\to 0^+} \br[\bigg]{ \frac{z^{-s+m-1}}{-s+m-1} } \br[\bigg]{ \frac{ 1 - \exp\pr[\big]{ z \logNorm{\lap} } }{ -\logNorm{\lap} } } \bigg\rvert_{z=w}^1 - \int_0^1 \br[\bigg]{ \frac{z^{-s+m-1}}{-s+m-1} } \exp(z\logNorm{\lap} ) \, dz \nonumber \\
& \qquad + \lim_{w\to\infty} \frac{z^{-s+m-1}}{-s+m-1} \bigg\rvert_{z=1}^w \\
& = \br[\bigg]{ \frac{1}{-s+m-1} } \br[\bigg]{ \frac{ 1 + \logNorm{\lap} - \exp\pr[\big]{ \logNorm{\lap} } }{ -\logNorm{\lap} } }
- \int_0^1 \br[\bigg]{ \frac{z^{-s+m-1}}{-s+m-1} } \exp(z\logNorm{\lap} ) \, dz \nonumber \\
& \le \br[\bigg]{ \frac{1}{-s+m-1} } \br[\bigg]{ \frac{ 1 + \logNorm{\lap} - \exp\pr[\big]{ \logNorm{\lap} } }{ -\logNorm{\lap} } }
+ \int_0^1 \br[\bigg]{ \frac{z^{-s+m-1}}{s-m+1} } \, dz \nonumber \\
& = \br[\bigg]{ \frac{1}{-s+m-1} } \br[\bigg]{ \frac{ 1 + \logNorm{\lap} - \exp\pr[\big]{ \logNorm{\lap} } }{ -\logNorm{\lap} } } + \frac{1}{(-s+m)(s-m+1)} < \infty . \nonumber
\end{align}
Combining \cref{eq:2.18_int}, \cref{eq:2.19_int}, 
and the fact that \cref{prop:high_lap_defined} assures that for all $u \in \ellTwo$ it holds that $\normm{ (-\lap)^m u } < \infty$
hence proves \cref{eq:2.14_int}.
\end{aproof}

\cfclear
\begin{definition}[Discrete fractional Laplace operator for $s\in(0,1)$]\label{def:frac_lap_pre}
\cfconsiderloaded{def:frac_lap_pre}
\cfload
Let $s \in (0,1)$.
Then we denote by $\preFrac \colon \allowbreak \ellTwo \allowbreak \to \ellTwo$ 
the function which satisfies for all $u \in \ellTwo$, $n\in\Z$ that
\begin{equation}
\pr[\big]{ \preFrac u }(n) = \frac{1}{\gamfn(-s)} \int_0^\infty z^{-s-1} \br[\big]{ \semigroup_z(\lap) - \id } u(n) \, dz \ifnocf.
\end{equation}
\cfout[.]
\end{definition}

\cfclear
\begin{definition}[Discrete fractional Laplace operator for $s\in(0,\infty)$]\label{def:frac_lap}
\cfconsiderloaded{def:frac_lap}
\cfload
Let $s \in (0,\infty)$.
Then we denote
by $\fracLap \colon \allowbreak \ellTwo \allowbreak \to \ellTwo$ the function which satisfies 
for all $u \in \ellTwo$, $n\in\Z$ that
\begin{align}\label{eq:def_frac_lap}
\pr[\big]{ \fracLap u }(n) & = \pr[\big]{ (-\Delta)^{s-\floor{s}} (-\Delta)^{\floor{s}} u}(n) \\
& = 
\begin{cases}
\pr[\big]{ \preFrac u}(n) & \colon s \in \N \\
\frac{1}{\gamfn(-(s-\floor{s}))} \int_0^\infty z^{-(s-\floor{s})-1} \br[\big]{ \semigroup_z(\lap) - \id } \pr[\big]{(-\Delta)^{\floor{s}} u}(n) \, 
dz & \colon s \in (0,\infty) \backslash \N
\end{cases} \ifnocf. \nonumber
\end{align}
\cfout[.]
\end{definition}


\section{The discrete fractional kernel}\label{sec:4}

In this section we introduce a kernel function which will allow us to conveniently provide a series representation of \cref{eq:def_frac_lap} in \cref{def:frac_lap}.
\cref{lem:gam_pre,lem:term_bd} are preliminary results which allow us to prove \cref{lem:kern_bd}---a result which outlines useful properties exhibited by the kernel defined in \cref{def:kern}.
It is worth noting that \cref{lem:gam_pre} is a well-known result and that \cref{lem:term_bd} is a generalization of
Ciaurri et al.\ \cite[Lemma 9.2 (a)]{ciaurri2018nonlocal}, which was only proven in the case where $s\in(0,1)$.
\cref{prop:gamma,lem:sum} are the main results of this section and allow us to prove \cref{lem:rep1} in \cref{sec:5_1}.

\cfclear
\begin{definition}[Fractional kernel]\label{def:kern}
\cfconsiderloaded{def:kern}
\cfload
We denote\footnote{Let $A \subseteq \R$. Then it holds for all $x \in A$ that $\1_A(x) = 1$ and for all $x \in \R \backslash A$ that $\1_A(x) = 0$.}
by $\Kern_s \colon \Z \to \R$, $s \in \R$, the function which satisfies for all $k \in \Z$, $m \in \N$, $s \in (m-1,m)$ that
\begin{equation}
\Kern_s(k) = \frac{ - \1_{\Z\backslash\{0\}}(k) \, 4^s \gamfn(\nicefrac{1}{2} + s) \gamfn( \abs{k} -s ) }{ \sqrt{\pi} 
\gamfn(-s) \gamfn( \abs{k} + 1 + s ) } \ifnocf.
\end{equation}
\cfout[.]
\end{definition}

\cfclear
\begin{athm}{proposition}{lem:gam_pre}
It holds for all $a,b\in\R$, $\lambda \in (0,\infty)$ with $0 \le a < b < \infty$ that
\begin{equation}\label{eq:3.2}
\min\{ \lambda , 1 \} \le \frac{ b^\lambda - a^\lambda }{ b^{\lambda-1} ( b-a ) } \le \max\{ \lambda, 1\}.
\end{equation}
\end{athm}

\begin{aproof}
First, note that for all $a,b\in\R$ with $0 \le a < b < \infty$ it holds that 
\begin{equation}\label{eq:3.3}
0 \le \nicefrac{a}{b} < 1
\end{equation}
Observe that \cref{eq:3.3} ensures that for all $\lambda \in (0,1)$, $a,b\in\R$ with $0 \le a < b < \infty$ it holds that
$0 \le \nicefrac{a}{b} \le \pr[]{ \nicefrac{a}{b} }^\lambda < 1$.
This assures that for all $ \lambda \in (0,1)$, $a,b\in\R$ with $0 \le a < b < \infty$ it holds that
\begin{equation}\label{eq:3.4}
\frac{ b^\lambda - a^\lambda }{ b^{\lambda-1} ( b-a ) } = \frac{ 1 - \pr[]{ \nicefrac{a}{b} }^\lambda }{ 1 - \nicefrac{a}{b} }  \le 1.
\end{equation}
Combining the mean value theorem and \cref{eq:3.4} hence guarantee that for all $ \lambda \in (0,1)$, $a,b\in\R$ with $0 \le a < b < \infty$
it holds that
there exists $c \in (\nicefrac{a}{b},1)$ such that
\begin{equation}\label{eq:3.5}
\frac{ b^\lambda - a^\lambda }{ b^{\lambda-1} ( b-a ) } = \frac{ 1 - \pr[]{ \nicefrac{a}{b} }^\lambda }{ 1 - \nicefrac{a}{b} } = \lambda c^{\lambda-1} \ge \lambda.
\end{equation}
Next, note that \cref{eq:3.3} demonstrates that for all $\lambda \in [1,\infty)$, $a,b\in\R$ with $0 \le a < b < \infty$ it holds that
$0 \le \pr[]{ \nicefrac{a}{b} }^\lambda \le \nicefrac{a}{b} < 1$.
This shows that for all $ \lambda \in [1,\infty)$, $a,b\in\R$ with $0 \le a < b < \infty$ it holds that
\begin{equation}\label{eq:3.6}
\frac{ b^\lambda - a^\lambda }{ b^{\lambda-1} ( b-a ) } = \frac{ 1 - \pr[]{ \nicefrac{a}{b} }^\lambda }{ 1 - \nicefrac{a}{b} }  \ge 1.
\end{equation}
Combining \cref{eq:3.6} with the mean value theorem therefore proves that 
for all $ \lambda \in [1,\infty)$, $a,b\in\R$ with $0 \le a < b < \infty$
it holds that
there exists $d \in (\nicefrac{a}{b},1)$ such that
\begin{equation}\label{eq:3.7}
\frac{ b^\lambda - a^\lambda }{ b^{\lambda-1} ( b-a ) } = \frac{ 1 - \pr[]{ \nicefrac{a}{b} }^\lambda }{ 1 - \nicefrac{a}{b} } = \lambda d^{\lambda-1} \le \lambda.
\end{equation}
Combining \cref{eq:3.4,eq:3.5,eq:3.6,eq:3.7} thus establishes \cref{eq:3.2}.
\end{aproof}

\cfclear
\begin{athm}{lemma}{lem:term_bd}
It holds for all $ m \in \N$, $s\in(m-1,m)$ that there exists $C \in \R$ such that for all $k \in \Z$ with $\abs{k} \in [m , \infty)$ it holds that
\begin{equation}\label{eq:3.8}
\abs[\Bigg]{ \frac{\gamfn(\abs{k} - s) }{\gamfn( \abs{k} + 1 + s)} - \frac{1}{\abs{k}^{1+2s} } } \le \frac{ C }{ \abs{k}^{2+2s} } \ifnocf.
\end{equation}
\cfload.
\end{athm}

\begin{aproof}
Throughout this proof let $m \in \N$, $s\in(m-1,m)$
and without loss of generality let $k \in \Z$ with $k \in [m,\infty)$.
Note that the triangle inequality assures that 
\begin{align}\label{eq:3.9}
\abs[\Bigg]{ \frac{\gamfn(\abs{k} - s) }{\gamfn( \abs{k} + 1 + s)} - \frac{1}{\abs{k}^{1+2s} } }
& \le 
\abs[\Bigg]{ \frac{\gamfn(\abs{k} - s) }{\gamfn( \abs{k} + 1 + s)} - \frac{1}{\abs{k-s}^{1+2s} } } 
+ \abs[\Bigg]{ \frac{1}{\abs{k-s}^{1+2s} } - \frac{1}{\abs{k}^{1+2s} } } \\
& = \abs[\Bigg]{ \frac{\gamfn(k - s) }{\gamfn( k + 1 + s)} - \frac{1}{\pr{k-s}^{1+2s} } } 
+ \abs[\Bigg]{ \frac{1}{\pr{k-s}^{1+2s} } - \frac{1}{k^{1+2s} } } . \nonumber
\end{align}
Next, observe that \cref{lem:gam_pre} (applied 
with $\lambda \with 1 + 2s$, $a \with \nicefrac{1}{k}$, $b \with \nicefrac{1}{(m-s)}$ in the notation of \cref{lem:gam_pre}) ensures that 
\begin{align}
& \abs[\Bigg]{ \frac{1}{\pr{k-s}^{1+2s} } - \frac{1}{k^{1+2s} } }
= \br[\Bigg]{ \frac{\abs[\big]{ \pr{k-s}^{-(1+2s)} - k^{-(1+2s)} } }{ \pr{k-s}^{-2s} \br[\big]{ \pr{k-s}^{-1} - k^{-1} } } } \pr{k-s}^{-2s} \br[\big]{ \pr{k-s}^{-1} - k^{-1} } \nonumber \\
& \quad \le \frac{\max\{ 1 + 2s, 1 \} }{\pr{k-s}^{2s}} \br[\Bigg]{ \frac{1}{k-s} - \frac{1}{k} }
= \frac{1 + 2s}{\pr{k-s}^{2s}} \br[\Bigg]{ \frac{k - \pr{k-s}}{\pr{k-s} k } }
= \frac{(1+2s)s}{k \pr{ k-s }^{1+2s} }.
\end{align}
This and the fact that $k \in [m,\infty)$ show that
\begin{equation}\label{eq:3.11}
\abs[\Bigg]{ \frac{1}{\pr{k-s}^{1+2s} } - \frac{1}{k^{1+2s} } }
\le \frac{(1+2s)s}{ \pr{ k-s }^{2+2s} }
\le \left[ \sup_{ k \in [m,\infty) } \frac{ k^{2+2s} }{ (k-s)^{2+2s} } \right] \frac{ (1+2s) s }{ k^{2+2s} }.
\end{equation}
In addition, note that, e.g., Tricomi and Erd{\'e}lyi \cite[Eq.\ (15), page 140]{tricomi1951asymptotic} (applied with $z \with k$, $\alpha \with -s$, $\beta \with 1+s$ in the notation of Tricomi and Erd{\'e}lyi \cite[Eq.\ (15), page 140]{tricomi1951asymptotic}) guarantees that
\begin{equation}
\frac{ \gamfn( k - s) }{\gamfn(k + 1 + s)} = \frac{1}{\gamfn( 1 + 2s )} \int_0^\infty \exp( -(k-s) v ) \pr[\big]{ 1 - \exp(-v) }^{2s} \, dv .
\end{equation}
This, the fact that for all $r \in [0,\infty)$ it holds that
$\int_0^\infty \exp( - r v ) v^{2s} \, dv = \gamfn(1+2s) r^{-(1+2s)}$,
and Jensen's inequality prove that
\begin{equation}\label{eq:3.13}
\begin{split}
& \abs[\Bigg]{ \frac{\gamfn(k - s) }{\gamfn( k + 1 + s)} - \frac{1}{\pr{k-s}^{1+2s} } } \\
& \quad = \abs[\Bigg]{ \frac{1}{\gamfn( 1 + 2s )} \int_0^\infty \exp( -(k-s) v ) \pr[\big]{ 1 - \exp(-v) }^{2s} \, dv - \frac{1}{\pr{k-s}^{1+2s} } } \\
& \quad = \abs[\Bigg]{ \frac{1}{\gamfn( 1 + 2s )} \int_0^\infty \exp( -(k-s) v ) \br[\Big]{ \pr[\big]{ 1 - \exp(-v) }^{2s} - v^{2s} } \, dv } \\
& \quad \le \frac{1}{\gamfn( 1 + 2s )} \int_0^\infty \exp( -(k-s) v ) v^{2s} \abs[\bigg]{ \pr[\Big]{ \frac{1 - \exp(-v)}{v} }^{\!2s} - 1 } \, dv . 
\end{split}
\end{equation}
Moreover, observe that \cref{lem:gam_pre} (applied with $\lambda \with 2s$, $a \with \nicefrac{(1-\exp(-v))}{v}$, $b \with 1$ in the notation of \cref{lem:gam_pre}) ensures that for all $v \in (0,\infty)$ it holds that
\begin{equation}
1 - \pr[\Big]{ \frac{1 - \exp(-v)}{v} }^{\!2s} \le \max\{2s,1\} \br[\bigg]{ 1 - \frac{1-\exp(-v)}{v} }.
\end{equation}
Combining this, \cref{eq:3.13}, the fact that for all $v\in(0,\infty)$ it holds that $v - (1-\exp(-v)) < \nicefrac{v^2}{2}$,
the fact that for all $r \in [0,\infty)$ it holds that
$\int_0^\infty \exp( - r v ) v^{1+2s} \, dv = \gamfn(2+2s) r^{-(2+2s)}$,
and the fact that \cref{def:gamma} implies that for all $z\in(0,\infty)$ it holds that $\gamfn(1+z) = z\gamfn(z)$ 
assures that
\begin{equation}
\begin{split}
& \abs[\Bigg]{ \frac{\gamfn(k - s) }{\gamfn( k + 1 + s)} - \frac{1}{\pr{k-s}^{1+2s} } } \\
& \quad \le \frac{\max\{2s,1\}}{\gamfn( 1 + 2s )} \int_0^\infty \exp( -(k-s) v ) v^{2s} \abs[\bigg]{ 1-  \pr[\Big]{ \frac{1 - \exp(-v)}{v} } } \, dv \\
& \quad \le \frac{\max\{s,\nicefrac{1}{2}\}}{\gamfn( 1 + 2s )} \int_0^\infty \exp( -(k-s) v ) v^{1+2s} \, dv
= \frac{\max\{s,\nicefrac{1}{2}\} \gamfn(2+2s)}{\gamfn(1+2s)} \frac{1}{\pr{k-s}^{2+2s}}  \\
& \quad = \max\{s,\nicefrac{1}{2}\}(1+2s) \frac{1}{\pr{k-s}^{2+2s}} \le \left[ \sup_{k\in[m,\infty)} \frac{k^{2+2s}}{\pr{k-s}^{2+2s}} \right] \frac{ \max\{s,\nicefrac{1}{2}\}(1+2s) }{ k^{2+2s} } .
\end{split}
\end{equation}
Combining this, \cref{eq:3.9}, \cref{eq:3.11}, and the fact that $\sup_{k\in[m,\infty)} k^{2+2s} (k-s)^{-(2+2s)} \in \R$ hence proves \cref{eq:3.8}.
\end{aproof}

\cfclear
\begin{athm}{lemma}{lem:kern_bd}
Let $ m \in \N$, $s\in(m-1,m)$.
Then
\begin{enumerate}[label=(\roman*)]
\item\label{lem:kern_bd_i1} it holds for all $k \in \Z$ that $\Kern_s(-k) = \Kern_s(k)$ and
\item\label{lem:kern_bd_i2} it holds that there exists $C \in \R$ such that for all $k \in \Z$ with $\abs{k} \in \Z \backslash \{0\}$ it holds that
\begin{equation}\label{eq:3.16}
\abs{ \Kern_s(k) } \le \frac{C}{\abs{k}^{1+2s}} \ifnocf.
\end{equation}
\end{enumerate}
\cfout.
\end{athm}

\begin{aproof}
First, note that for all $k \in \Z$ it holds that
\begin{align}
\Kern_s(-k) & = \frac{ - \1_{\Z\backslash\{0\}}(-k) \, 4^s \gamfn(\nicefrac{1}{2} + s) \gamfn( \abs{-k} -s ) }{ \sqrt{\pi} 
\gamfn(-s) \gamfn( \abs{-k} + 1 + s ) } \\
& = \frac{ - \1_{\Z\backslash\{0\}}(k) \, 4^s \gamfn(\nicefrac{1}{2} + s) \gamfn( \abs{k} -s ) }{ \sqrt{\pi} 
\gamfn(-s) \gamfn( \abs{k} + 1 + s ) } = \Kern_s(k) \nonumber
\end{align}
\cfload.
This establishes \cref{lem:kern_bd_i1}.
Next, observe that for all $k \in \Z \cap (-m,m)$ with $k \neq 0$ it holds that
\begin{align}
\abs{ \Kern_s(k) } & \le \frac{ 4^s \gamfn(\nicefrac{1}{2} + s) \abs{ \gamfn( \abs{k} -s ) } }{ \sqrt{\pi} 
\abs{ \gamfn(-s) } \gamfn( \abs{k} + 1 + s ) } 
= \frac{ 4^s \gamfn(\nicefrac{1}{2} + s) \abs{ \gamfn( \abs{k} -s ) } }{ \sqrt{\pi} 
\abs{ \gamfn(-s) } \gamfn( \abs{k} + 1 + s ) } \cdot \frac{\abs{k}^{1+2s}}{\abs{k}^{1+2s}} \nonumber \\
& \le \frac{ 4^s \gamfn( \nicefrac{1}{2} + s) m^{1+2s}}{ \sqrt{\pi} \abs{\gamfn(-s)} } \br[\Bigg]{ \sup_{j \in \Z \cap (-m,m) } \frac{ \abs{\gamfn( \abs{j} - s )} }{ \gamfn( \abs{j} + 1 + s )} } \frac{1}{ \abs{k}^{1+2s}} .
\end{align}
This, the fact that $s \in (0,\infty) \backslash \N$ implies that $\sup_{j \in \Z \cap (-m,m) } \nicefrac{ \abs{\gamfn( \abs{j} - s )} }{ \gamfn( \abs{j} + 1 + s )} \in \R$, and \cref{lem:term_bd} establish \cref{lem:kern_bd_i2}.
\end{aproof}

\cfclear
\begin{athm}{proposition}{prop:gamma}
It holds for all $m \in \N$, $s\in(m-1,m)$, $k \in \Z$ that
\begin{equation}\label{eq:gamma_prop}
\Kern_s(k) = \frac{ \1_{\Z\backslash\{0\}}(k) \, (-1)^{k+1} \gamfn(2s+1) }{ \gamfn(1+s+k) \gamfn(1+s-k) }
\end{equation}
\cfload.
\end{athm}

\begin{aproof}
Throughout this proof let $m \in \N$, $s\in(m-1,m)$ and without loss of generality let $k \in \N$ (cf.\ \cref{lem:kern_bd_i1} of \cref{lem:kern_bd}).
Observe that the Legendre duplication formula 
(cf., e.g., Abramowitz and Stegun \cite[Eq.\ (6.1.18), Page 256]{abramowitz1970handbook}) ensures that
\begin{align}
\Kern_s(k) & = \frac{ - 4^s \gamfn(\nicefrac{1}{2} + s) \gamfn( k -s ) }{ \sqrt{\pi} 
\gamfn(-s) \gamfn( k + 1 + s ) }
= \frac{ - 4^s \gamfn(\nicefrac{1}{2} + s) \gamfn( k -s ) }{ \sqrt{\pi} 
\gamfn(-s) \gamfn( k + 1 + s ) } \cdot \frac{\gamfn(s)}{\gamfn(s)} \nonumber \\
& = \frac{ - 4^s \br[\big]{ 2^{1-2s} \sqrt{\pi} \gamfn(2s) } \gamfn( k -s ) }{ \sqrt{\pi} 
\gamfn(-s) \gamfn(s) \gamfn( k + 1 + s ) }
= \frac{ -2 \gamfn(2s) \gamfn(k-s) }{ \gamfn(-s) \gamfn(s) \gamfn( k + 1 + s ) }.
\end{align}
This and the fact that \cref{def:gamma} implies that for all $z\in \C$ with $\real(z) \in \R \backslash \{\ldots, -2,-1,0\}$ it holds that $z\gamfn(z) = \gamfn(1+z)$ assure that
\begin{equation}\label{eq:gamma_prop1}
\Kern_s(k)
= \frac{ -2 \gamfn(2s) \gamfn(k-s) }{ \gamfn(-s) \gamfn(s) \gamfn( k + 1 + s ) } \cdot \frac{s}{s} 
= \frac{ \gamfn(2s+1) \gamfn(k-s) }{ \gamfn(1-s) \gamfn(s) \gamfn( k + 1 + s ) }.
\end{equation}
Next, note that the Euler reflection formula (cf., e.g., Abramowitz and Stegun \cite[Eq.\ (6.1.17), Page 256]{abramowitz1970handbook}) guarantees that
\begin{equation}\label{eq:gamma_prop2}
\gamfn(s) \gamfn(1-s) = (-1)^{k+1} \gamfn(k-s) \gamfn(1+s-k) .
\end{equation}
Combining \cref{eq:gamma_prop1,eq:gamma_prop2}
hence yields \cref{eq:gamma_prop}.
\end{aproof}

\cfclear
\begin{athm}{proposition}{prop:sum}
Let $s \in (0,\infty) \backslash \N$.
Then it holds for all $m\in\N$ that
\begin{equation}\label{eq:sum_claim}
\frac{\gamfn(m-s)}{2s\gamfn(m+s)} + \sum_{k=1}^{m-1} \frac{ \gamfn( k - s) }{ \gamfn( k + 1 + s )}
= \frac{ - \gamfn(-s) }{ 2 \gamfn(1+s) }
\end{equation}
\cfout.
\end{athm}

\begin{aproof}
We prove \cref{eq:sum_claim} by induction on $m \in \N$. 
For the base case $m=1$ note that the fact that \cref{def:gamma} ensures that for all $z \in \C$ with $\real(z) \in \R \backslash\{\ldots,-2,-1,0\}$ it holds that $z\gamfn(z) = \gamfn(z+1)$ guarantees that 
\begin{equation}
\frac{\gamfn(1-s)}{2s\gamfn(1+s)} + \sum_{k=1}^{0} \frac{ \gamfn( k - s) }{ \gamfn( k + 1 + s )} 
= \frac{\gamfn(1-s)}{2s\gamfn(1+s)} 
= \frac{-s\gamfn(-s)}{2s\gamfn(1+s)}
= \frac{ -\gamfn(-s) }{ 2 \gamfn(1+s) }
\end{equation}
\cfload.
This establishes \cref{eq:sum_claim} in the case $m=1$. For the induction step $ \N \ni (m-1) \induct m \in \N \cap [2,\infty)$, let $m \in \N \cap [2,\infty)$ and assume for all $\mathfrak{m} \in \{1,2,\ldots,m-1\}$
that \cref{eq:sum_claim} holds. Observe that the induction hypothesis 
shows that for all $m \in \N \cap [2,\infty)$
it holds that
\begin{align}\label{eq:sum_claim1}
& \frac{\gamfn(m-s)}{2s\gamfn(m+s)} + \sum_{k=1}^{m-1} \frac{ \gamfn( k - s) }{ \gamfn( k + 1 + s )} \nonumber \\
& \quad = \left[ \frac{\gamfn(m-s)}{2s\gamfn(m+s)} + \frac{\gamfn((m-1)-s)}{\gamfn((m-1)+1+s)} \right] + \sum_{k=1}^{(m-1)-1} \frac{ \gamfn( k - s) }{ \gamfn( k + 1 + s )} \\
& \quad = \frac{\gamfn(m-s)}{2s\gamfn(m+s)} + \frac{\gamfn(m-1-s)}{\gamfn(m+s)}
+ \frac{-\gamfn(m-1-s)}{2s\gamfn(m-1+s)} + \frac{-\gamfn(-s)}{2\gamfn(1+s)} . \nonumber
\end{align}
Next, note that the fact that \cref{def:gamma} ensures that for all $z \in \C$ with $\real(z) \in \R \backslash\{\ldots,-2,\allowbreak -1, \allowbreak 0\}$ it holds that $z\gamfn(z) = \gamfn(z+1)$  demonstrates that for all $m \in \N \cap [2,\infty)$ it holds that
\begin{align}\label{eq:sum_claim2}
& \frac{\gamfn(m-s)}{2s\gamfn(m+s)} + \frac{-\gamfn(m-1-s)}{2s\gamfn(m-1+s)}
= \frac{\gamfn(m-s)}{2s\gamfn(m+s)} + \frac{-(m-1+s)\gamfn(m-1-s)}{2s(m-1+s)\gamfn(m-1+s)} \nonumber \\
& \quad = \frac{ \gamfn(m-s) - (m-1+s) \gamfn(m-1-s)}{2s\gamfn(m+s)}
= \frac{ (m-1-s)\gamfn(m-s) - (m-1+s) \gamfn(m-s)}{2s(m-1-s)\gamfn(m+s)} \nonumber \\
& \quad = \left[ \frac{(m-1-s) - (m-1+s)}{2s(m-1-s)} \right] \frac{\gamfn(m-s)}{\gamfn(m+s)}
= \frac{-\gamfn(m-s)}{(m-1-s)\gamfn(m+s)}.
\end{align}
Moreover, observe that the fact that
\cref{def:gamma} assures that for all $z \in \C$ with $\real(z) \in \R \backslash\{\ldots,-2,\allowbreak -1, \allowbreak 0\}$ it holds that $z\gamfn(z) = \gamfn(z+1)$  demonstrates that for all $m \in \N \cap [2,\infty)$ it holds that
\begin{align}\label{eq:sum_claim3}
\frac{-\gamfn(m-s)}{(m-1-s)\gamfn(m+s)} + \frac{\gamfn(m-1-s)}{\gamfn(m+s)}
& = \frac{ -\gamfn(m-s) + (m-1-s)\gamfn(m-1-s)}{(m-1-s)\gamfn(m+s)} \nonumber \\
& = \frac{ -\gamfn(m-s) + \gamfn(m-s)}{(m-1-s)\gamfn(m+s)} = 0.
\end{align}
Combining \cref{eq:sum_claim1,eq:sum_claim2,eq:sum_claim3} therefore proves \cref{eq:sum_claim}.
\end{aproof}

\cfclear
\begin{athm}{lemma}{lem:sum}
It holds for all $m\in\N$, $s\in(m-1,m)$ that
\begin{equation}\label{eq:sum}
\sum_{k\in\Z} \Kern_s(k) 
= \frac{ 4^s \gamfn(\nicefrac{1}{2}+s)}{\sqrt{\pi}\gamfn(1+s)}
\end{equation}
\cfload.
\end{athm}

\begin{aproof}
First, note that \cref{lem:kern_bd_i2} of \cref{lem:kern_bd} ensures that for all $m \in\N$, $s\in(m-1,m)$ it holds that $\sum_{k\in\Z} \Kern_s(k) \in \R$.
Next, observe that \cref{lem:kern_bd_i1} of \cref{lem:kern_bd} 
assures that for all $m \in \N$, $s \in (m-1,m)$ it holds that
\begin{align}\label{eq:3_19}
\sum_{k\in\Z} \Kern_s(k) & = \frac{- 4^s \gamfn( \nicefrac{1}{2} + s)}{ \sqrt{\pi} \gamfn(-s) } \left[ \sum_{k\in\Z} \frac{ \1_{\Z\backslash\{0\}}(k) \gamfn( \abs{k} - s) }{ \gamfn( \abs{k} + 1 + s )} \right] \\
& = \frac{- 2\cdot 4^s \gamfn( \nicefrac{1}{2} + s)}{ \sqrt{\pi} \gamfn(-s) } \left[ \sum_{k\in\N} \frac{ \gamfn( k - s) }{ \gamfn( k + 1 + s )} \right]. \nonumber
\end{align}
In addition, note that, e.g., Artin \cite[Eq.\ (2.13)]{artin2015gamma} (applied for every $m \in \N$, $s \in (m-1,m)$,
$k \in \N \cap [m,\infty)$ with $x \with k-s$, $y \with 1+2s$ in the notation of Artin \cite[Eq.\ (2.13)]{artin2015gamma}) implies that
for all 
$m \in \N$, $s \in (m-1,m)$,
$k \in \N \cap [m,\infty)$ it holds that
\begin{align}\label{eq:beta_fn}
\frac{ \gamfn( k - s) }{ \gamfn( k + 1 + s )} 
= \frac{1}{\gamfn(1+2s)} \br[\bigg]{ \frac{ \gamfn( k - s) \gamfn(1+2s) }{ \gamfn( k + 1 + s )}  }
& = \frac{1}{\gamfn( 1+2s )} \int_0^1 (1-z)^{(1+2s)-1} z^{(k-s)-1} \, dz \nonumber \\
& = \frac{1}{\gamfn(1+2s)} \int_0^1 (1-z)^{2s} z^{k-s-1} \, dz.
\end{align}
This and Fubini's theorem
guarantee that for all $m \in \N$, $s \in (m-1,m)$
it holds that
\begin{align}\label{eq:3.21_gamma}
& \sum_{k\in\N} \frac{ \gamfn( k - s) }{ \gamfn( k + 1 + s )}
= \sum_{k=1}^{m-1} \frac{ \gamfn( k - s) }{ \gamfn( k + 1 + s )} + \sum_{k=m}^\infty \frac{ \gamfn( k - s) }{ \gamfn( k + 1 + s )} \nonumber \\
& \quad = \sum_{k=1}^{m-1} \frac{ \gamfn( k - s) }{ \gamfn( k + 1 + s )} + \sum_{k=m}^\infty \left[ \frac{1}{\gamfn(1+2s)} \int_0^1 (1-z)^{2s} z^{k-s-1} \, dz \right] \nonumber \\
& \quad = \sum_{k=1}^{m-1} \frac{ \gamfn( k - s) }{ \gamfn( k + 1 + s )} + \frac{1}{\gamfn(1+2s)} \int_0^1 (1-z)^{2s} z^{-s-1} \left[ \sum_{k=m}^\infty z^{k} \right] \, dz \\
& \quad = \sum_{k=1}^{m-1} \frac{ \gamfn( k - s) }{ \gamfn( k + 1 + s )} + \frac{1}{\gamfn(1+2s)} \int_0^1 (1-z)^{2s} z^{m-s-1} \left[ \sum_{k=0}^\infty z^{k} \right] \, dz \nonumber \\
& \quad = \sum_{k=1}^{m-1} \frac{ \gamfn( k - s) }{ \gamfn( k + 1 + s )} + \frac{1}{\gamfn(1+2s)} \int_0^1 (1-z)^{2s-1} z^{m-s-1} \, dz. \nonumber 
\end{align}
Combining this, the fact that \cref{def:gamma} ensures that for all $z \in \C$ with $\real(z) \in \R \backslash\{\ldots,-2,-1,0\}$ it holds that $z\gamfn(z) = \gamfn(z+1)$, and, e.g., Artin \cite[Eq.\ (2.13)]{artin2015gamma} (applied for every $m \in \N$, $s \in (m-1,m)$ with $x \with m-s$, $y \with 2s$ in the notation of Artin \cite[Eq.\ (2.13)]{artin2015gamma})
demonstrates that for all $m \in \N$, $s \in (m-1,m)$
it holds that
\begin{align}
\sum_{k\in\N} \frac{ \gamfn( k - s) }{ \gamfn( k + 1 + s )}
& = \sum_{k=1}^{m-1} \frac{ \gamfn( k - s) }{ \gamfn( k + 1 + s )} + \frac{1}{\gamfn(1+2s)} \left[ \frac{\gamfn(m-s)\gamfn(2s)}{\gamfn(m+s)} \right]\nonumber \\
& = \sum_{k=1}^{m-1} \frac{ \gamfn( k - s) }{ \gamfn( k + 1 + s )} + \frac{1}{2s\gamfn(2s)} \left[ \frac{\gamfn(m-s)\gamfn(2s)}{\gamfn(m+s)} \right] \\
& = \sum_{k=1}^{m-1} \frac{ \gamfn( k - s) }{ \gamfn( k + 1 + s )} + \frac{\gamfn(m-s)}{2s\gamfn(m+s)}. \nonumber
\end{align}
Combining
this and \cref{prop:sum} proves that for all $m \in \N$, $s\in (m-1,m)$ it holds that
\begin{align}
\sum_{k\in\Z} \Kern_s(k) & = \frac{- 2\cdot 4^s \gamfn( \nicefrac{1}{2} + s)}{ \sqrt{\pi} \gamfn(-s) } \left[ \sum_{k=1}^{m-1} \frac{ \gamfn( k - s) }{ \gamfn( k + 1 + s )} + \frac{\gamfn(m-s)}{2s\gamfn(m+s)} \right] \\
& = \frac{ - 2\cdot 4^s \gamfn( \nicefrac{1}{2} + s)}{ \sqrt{\pi} \gamfn(-s) } \left[ \frac{-\gamfn(-s)}{2\gamfn(1+s)} \right] = \frac{4^s \gamfn( \nicefrac{1}{2} + s )}{ \sqrt{\pi} \gamfn(1+s) }. \nonumber
\end{align}
\end{aproof}


\section{Series representation for the discrete fractional Laplace operator}\label{sec:5}

In this section we prove the main result of the article.
First, in \cref{sec:5_1} we provide a series representation of the real-valued non-integer powers of the discrete Laplace operator (cf.\ \cref{lem:rep1}).
This representation employs the fractional kernel introduced in \cref{def:kern} and its proof hinges upon the results developed in \cref{sec:4}.
It is particularly interesting to note that the representation obtained in \cref{lem:rep1} coincides with the representation presented in Ciaurri et al.\ \cite[Theorem 1.1]{ciaurri2018nonlocal} (the case where $m=1$) and Padgett et al.\ \cite[Theorem 2]{padgett2020anomalous} (the case where $m=2$).

In \cref{sec:5_2} we demonstrate that the representation presented in \cref{lem:rep1} holds for $s\in\N$ if we consider the limiting values of the discrete kernel function.
The main result of the article, \cref{cor:final} below, follows immediately from the combination of \cref{lem:rep1,lem:ver2}.
In particular, \cref{cor:final} demonstrates that all real-valued positive powers of the discrete Laplace operator may be represented with the same series (or, at least, as the limit of this series).
Therefore, it is the case that the discrete fractional Laplace operator is, in some sense, a perturbation of the standard positive integer power case.

While it is not the purpose of this article to discuss such issues, we wish to emphasize the importance of the last sentence in the previous paragraph.
The fact that the discrete fractional Laplace operator's series representation coincides with the series representation for positive integer powers provides a framework to endow fractional calculus with potentially enlightening physical interpretations.
A particularly lacking feature of the fractional calculus is the lack of meaningful physical interpretations in many situations, which has been one of the primary limiting factor in its widespread application. 
However, \cref{cor:final} allows us to view the fractional powers as ``transitional phases'' between the positive integer cases, loosely speaking. 
Thus, we may use the known physical intuition for positive integer powers to provide a deeper understanding of the positive non-integer cases.

\subsection{Series representation for positive non-integer order}\label{sec:5_1}

\cfclear
\begin{athm}{lemma}{lem:rep1}
It holds for all $m \in \N$, $s\in(m-1,m)$, $u \in \ellTwo$, $n\in\Z$ that
\begin{equation}\label{eq:frac_rep}
\pr[\big]{ \fracLap u }(n) = \sum_{k\in\Z} \Kern_s(n-k) \pr[\big]{ u(n) - u(k) } 
\end{equation}
\cfout.
\end{athm}

\begin{aproof}
Throughout this proof let $u\in\ellTwo$, let $v \colon \Z \to \ellTwo$ satisfy for all $n\in\Z$ that 
$v(n) = (-\lap u)(n)$, and let $A_{s} \in \R$, $s\in(0,\infty)$,
satisfy for all $s\in(0,\infty)$ that
\begin{equation}
\label{eq:kern}
A_{s} = \frac{4^s \gamfn(\nicefrac{1}{2} + s)}{\sqrt{\pi} \gamfn(1 + s)} \ifnocf.
\end{equation}
\cfload.
Observe that Padgett et al.\ \cite[Theorem 1]{padgett2020anomalous} establishes \cref{eq:frac_rep} in the case that $m =1$, $s\in(0,1)$.
Next, note thatthe fact that for all $m \in \N \cap[2,\infty)$, $n\in\Z$ it holds that $((-\lap)^m u)(n) = ((-\lap)^{m-1}(-\lap) u)(n)$ (i.e., we are invoking the fact that standard function composition is associative on its domain) ensures that for all $m \in \N \cap[2,\infty)$, $s\in(m-1,m)$, $n\in\Z$ it holds that
\begin{align}\label{eq:5.3}
\pr[\big]{ \fracLap u}(n) & = \pr[\big]{ (-\Delta)^{s-(m-1)}\pr[]{ (-\lap)^{m-1} u } } (n) 
= \pr[\big]{ (-\Delta)^{s-(m-1)}\pr[]{ (-\lap)^{m-2} (-\lap u) } } (n) \nonumber \\
& = \pr[\big]{ (-\Delta)^{s-(m-1)}\pr[]{ (-\lap)^{m-2} v } } (n)
= \pr[\big]{ (-\lap)^{s-1} v }(n)
\end{align}
\cfload.
We now claim that for all $m\in\N\cap[2,\infty)$, $s\in(m-1,m)$, $n\in\N$ it holds that
\begin{equation}\label{eq:induct}
\pr[\big]{ (-\lap)^{s} u }(n)
= \sum_{k\in\Z} \Kern_s(n-k) \pr[\big]{ u(n) - u(k) }
\end{equation}
\cfload.
We prove \cref{eq:induct} by induction on $m \in \N \cap [2,\infty)$. 
For the base case $m=2$ observe that Padgett et al.\ \cite[Theorem 2]{padgett2020anomalous} establishes \cref{eq:induct}.
For the induction step $\N \cap [2,\infty) \ni (m-1) \induct m \in \N \cap [3,\infty)$, let $m \in \N \cap [3,\infty)$ and assume for all $\mathfrak{m} \in \{2,3,\dots,m-1\}$, $s \in (\mathfrak{m}-1,\mathfrak{m})$, $n\in\Z$ that \cref{eq:induct} holds.
Observe that the induction hypothesis, \cref{eq:kern}, \cref{eq:5.3}, \cref{eq:induct}, and \cref{lem:sum} 
demonstrate that for all $s\in(m-1,m)$, $n\in\Z$ it holds that
\begin{align}
& \pr[\big]{ (-\lap)^{s} u }(n)
= \pr[\big]{ (-\Delta)^{s-1} v }(n) 
= \SmallSum_{k\in\Z} \Kern_{s-1}(n-k) \pr[\big]{ v(n) - v(k) } \\
& \quad = v(n) \SmallSum_{k\in\Z} \Kern_{s-1}(n-k) - \SmallSum_{k\in\Z} \Kern_{s-1}(n-k) v(k) 
= A_{s-1} v(n) - \SmallSum_{k\in\Z} \Kern_{s-1}(n-k) v(k) . \nonumber
\end{align}
This, the fact that for all $n\in\Z$ it holds that $v(n) = (-\lap u)(n)$, and \cref{lem:kern_bd_i1} of \cref{lem:kern_bd} 
show that for all $s\in(m-1,m)$, $n\in\Z$ it holds that
\begin{align}\label{eq1}
\pr[\big]{ \fracLap u}(n) 
& = A_{s-1} v(n) - \SmallSum_{k\in\Z} \Kern_{s-1}(n-k) v(k) \nonumber \\
& = A_{s-1} \pr[\big]{ 2u(n) - u(n-1) - u(n+1) } \nonumber \\
& \qquad - \SmallSum_{k\in\Z} \Kern_{s-1}(n-k) \pr[\big]{ 2u(k) - u(k-1) - u(k+1) } \nonumber \\
& = A_{s-1} \pr[\big]{ 2u(n) - u(n-1) - u(n+1) } \nonumber \\
& \qquad - \SmallSum_{k\in\Z} \br[\big]{ 2\Kern_{s-1}(k) - \Kern_{s-1}(k-1) - \Kern_{s-1}(k+1) } u(n-k) \nonumber \\
& = A_{s-1} \pr[\big]{ 2u(n) - u(n-1) - u(n+1) } \nonumber \\
& \qquad - \br[\big]{ 2\Kern_{s-1}(0) - \Kern_{s-1}(-1) - \Kern_{s-1}(1) } u(n) \\
& \qquad - \br[\big]{ 2\Kern_{s-1}(1) - \Kern_{s-1}(0) - \Kern_{s-1}(2) } u(n-1) \nonumber \\
& \qquad - \br[\big]{ 2\Kern_{s-1}(-1) - \Kern_{s-1}(-2) - \Kern_{s-1}(0) } u(n+1) \nonumber \\
& \qquad - \SmallSum_{k\in\Z\backslash\{-1,0,1\}} \br[\big]{ 2\Kern_{s-1}(k) - \Kern_{s-1}(k-1) - \Kern_{s-1}(k+1) } u(n-k) \nonumber \\
& = \br[\big]{ 2 A_{s-1} + 2 \Kern_{s-1}(1) } u(n) \nonumber \\
& \qquad - \br[\big]{ A_{s-1} + 2\Kern_{s-1}(1) - \Kern_{s-1}(2) } \pr[\big]{ u(n-1) + u(n+1) } \nonumber \\
& \qquad - \SmallSum_{k\in\Z\backslash\{-1,0,1\}} \br[\big]{ 2\Kern_{s-1}(k) - \Kern_{s-1}(k-1) - \Kern_{s-1}(k+1) } u(n-k). \nonumber
\end{align}
Next, note that \cref{eq:kern}, the fact that \cref{def:gamma} implies that for all $z\in\C$ with $\real(z) \in \R \backslash\{\ldots,-2,\allowbreak -1, \allowbreak 0\}$ it holds that $\gamfn(z+1) = z\gamfn(z)$, and the fact that $1-s \in (-\infty,0]$ guarantee that for all $s\in(m-1,m)$ it holds that
\begin{align}\label{kern_exp1}
\Kern_{s-1}(1) & = \frac{ - \1_{\Z\backslash\{0\}}(1) \, 4^{s-1} \gamfn(\nicefrac{1}{2} + (s-1)) \gamfn( \abs{1} - (s-1) ) }{ \sqrt{\pi} \gamfn(-(s-1)) \gamfn( \abs{1} + 1 + (s-1) ) } 
= \frac{ - 4^{s-1} \gamfn(s - \nicefrac{1}{2} ) \gamfn( 2-s ) }{ \sqrt{\pi} \gamfn(1-s) \gamfn( 1 + s ) } \nonumber \\
& = \frac{ - 4^{s-1} \gamfn(s - \nicefrac{1}{2} ) (1-s) \gamfn( 1-s ) }{ \sqrt{\pi} \gamfn(1-s) s \gamfn( s ) } 
= \frac{ 4^{s-1} \gamfn(s-\nicefrac{1}{2}) (s-1) }{ \sqrt{\pi} s \gamfn( s ) } 
\end{align}
and
\begin{align}\label{kern_exp2}
\Kern_{s-1}(2) & = \frac{ - \1_{\Z\backslash\{0\}}(2) \, 4^{s-1} \gamfn(\nicefrac{1}{2} + (s-1)) \gamfn( \abs{2} - (s-1) ) }{ \sqrt{\pi} \gamfn(-(s-1)) \gamfn( \abs{2} + 1 + (s-1) ) } 
= \frac{ - 4^{s-1} \gamfn(s - \nicefrac{1}{2} ) \gamfn( 3-s ) }{ \sqrt{\pi} \gamfn(1-s) \gamfn( 2 + s ) } \nonumber \\
& = \frac{ - 4^{s-1} \gamfn(s - \nicefrac{1}{2} ) (2-s)(1-s) \gamfn( 1-s ) }{ \sqrt{\pi} \gamfn(1-s) (1+s) s \gamfn( s ) } 
= \frac{ 4^{s-1} \gamfn(s-\nicefrac{1}{2}) (2-s) (s-1) }{ \sqrt{\pi} (1+s) s \gamfn( s ) } .
\end{align}
Observe that \cref{eq:kern}, \cref{kern_exp1}, and the fact that \cref{def:gamma} implies that for all $z\in\C$ with $\real(z) \in \R \backslash\{\ldots,\allowbreak -2, \allowbreak -1, \allowbreak 0\}$ it holds that $\gamfn(z+1) = z\gamfn(z)$ hence ensure that for all $s\in(m-1,m)$ it holds that
\begin{align}\label{eq2}
2A_{s-1} + 2\Kern_{s-1}(1) & = 2 \br[\bigg]{ \frac{4^{s-1} \gamfn(s-\nicefrac{1}{2})}{\sqrt{\pi} \gamfn(s)} } + 2 \br[\bigg]{ \frac{ 4^{s-1} \gamfn(s-\nicefrac{1}{2}) (s-1) }{ \sqrt{\pi} s \gamfn( s ) } } \nonumber \\
& = \frac{ 2 \cdot 4^{s-1} \gamfn(s- \nicefrac{1}{2}) }{ \sqrt{\pi} \gamfn( s ) } \br[\bigg]{ 1 + \frac{ s-1 }{ s } } 
= \frac{ 2 \cdot 4^{s-1} \gamfn( s-\nicefrac{1}{2}) }{ \sqrt{\pi} \gamfn( s ) } \br[\bigg]{ \frac{ 2s-1 }{ s } } \\
& = \frac{ 4^{s} (s- \nicefrac{1}{2}) \gamfn(s-\nicefrac{1}{2})}{ \sqrt{\pi} s \gamfn(s) } 
= \frac{ 4^{s} \gamfn( \nicefrac{1}{2} + s ) }{ \sqrt{\pi} \gamfn( 1 + s ) } = A_{s}. \nonumber
\end{align}
In addition, note that \cref{prop:gamma} and the fact that \cref{def:gamma} implies that for all $z \in \C$ with $\real(z) \in \R \backslash \{\ldots,-2,-1,0\}$ it holds that $\gamfn(z+1) = z\gamfn(z)$ demonstrate that for all $k\in\Z \backslash \{-1,0,1\}$, $s\in(m-1,m)$ it holds that
\begin{align}\label{eq3}
& 2\Kern_{s-1}(k) - \Kern_{s-1}(k-1) - \Kern_{s-1}(k+1) \nonumber \\
& \quad = \frac{-2 (-1)^k \gamfn(2s-1)}{\gamfn(s+k)\gamfn(s-k) } - \frac{(s+k-1) (-1)^k \gamfn(2s-1)}{(s-k)\gamfn(s+k)\gamfn(s-k)} - \frac{(s-k-1) (-1)^k \gamfn(2s-1)}{(s+k)\gamfn(s+k)\gamfn(s-k)} \\
& \quad = \frac{(-1)^k \gamfn(2s-1)}{\gamfn(s+k)\gamfn(s-k)} \br[\bigg]{ -2 - \frac{s+k-1}{s-k} - \frac{s-k-1}{s+k} } \nonumber \\
& \quad = \frac{(-1)^k\gamfn(2s-1)}{\gamfn(s+k+1)\gamfn(s+k-1)} \br[\big]{ 2s-4s^2 }
= \frac{(-1)^{k+1} \gamfn(2s+1)}{\gamfn(1+s+k)\gamfn(1+s-k)} = \Kern_s(k). \nonumber 
\end{align}
Moreover, observe that \cref{eq:kern}, \cref{kern_exp1}, \cref{kern_exp2}, and the fact that \cref{def:gamma} implies that for all $z \in \C$ with $\real(z) \in \R \backslash \{\ldots,-2,-1,0\}$ it holds that $\gamfn(z+1) = z\gamfn(z)$ assure that for all $s \in (m-1,m)$ it holds that
\begin{align}\label{eq4}
& A_{s-1} + 2\Kern_{s-1}(1) - \Kern_{s-1}(2) \nonumber \\
& \quad = \frac{4^{s-1} \gamfn(s-\nicefrac{1}{2})}{\sqrt{\pi} \gamfn(s)} + 2 \br[\bigg]{ \frac{ 4^{s-1} \gamfn(s-\nicefrac{1}{2}) (s-1) }{ \sqrt{\pi} s \gamfn( s ) } } - \frac{ 4^{s-1} \gamfn(s-\nicefrac{1}{2}) (2-s) (s-1) }{ \sqrt{\pi} (1+s) s \gamfn( s ) } \nonumber \\
& \quad = \frac{4^{s-1} \gamfn(s-\nicefrac{1}{2})}{\sqrt{\pi} \gamfn(s)} \br[\Bigg]{ 1 + 2 \pr[\bigg]{ \frac{s-1}{s} } + \frac{(s-2)(s-1)}{s(1+s)} }
= \frac{4^{s-1} \gamfn(s-\nicefrac{1}{2})}{\sqrt{\pi} \gamfn(s)} \br[\Bigg]{ \frac{ 4 ( s-\nicefrac{1}{2} ) }{ s+1 } } \\
& \quad = \frac{ 4^{s} \gamfn( s+ \nicefrac{1}{2} )}{ \sqrt{\pi} \gamfn(s) (s+1) }
= \frac{ 4^{s} \gamfn( s+ \nicefrac{1}{2} )}{ \sqrt{\pi} \gamfn(s) (s+1) } \cdot \frac{ s\gamfn(-s) }{ s\gamfn(-s) }
= \frac{- 4^s \gamfn(s+\nicefrac{1}{2}) \gamfn(1-s)}{\sqrt{\pi} \gamfn(-s) \gamfn(2+s)} = \Kern_s(1). \nonumber
\end{align}
Combining \cref{eq1}, \cref{eq2}, \cref{eq3}, \cref{eq4}, and \cref{lem:kern_bd_i1} of \cref{lem:kern_bd} therefore yields that for all $s\in(m-1,m)$, $n\in\Z$ it holds that
\begin{align}
\pr[\big]{ \fracLap u}(n)
& = \br[\big]{ 2 A_{s-1} + 2 \Kern_{s-1}(1) } u(n) \nonumber \\
& \qquad - \br[\big]{ A_{s-1} + 2\Kern_{s-1}(1) - \Kern_{s-1}(2) } \pr[\big]{ u(n-1) + u(n+1) } \nonumber \\
& \qquad - \SmallSum_{k\in\Z\backslash\{-1,0,1\}} \br[\big]{ 2\Kern_{s-1}(k) - \Kern_{s-1}(k-1) - \Kern_{s-1}(k+1) } u(n-k) \\
& = A_s u(n) - \Kern_s(1) \pr[\big]{ u(n-1) + u(n+1) } -  \SmallSum_{k\in\Z\backslash\{-1,0,1\}} \Kern_s(k) u(n-k) \nonumber \\
& = A_s u(n) -  \SmallSum_{k\in\Z\backslash\{0\}} \Kern_s(k) u(n-k) 
= \SmallSum_{k\in\Z} \Kern_s(n-k) \pr[\big]{ u(n) - u(k) }. \nonumber
\end{align}
Induction hence establishes \cref{eq:induct}.
\end{aproof}

\subsection{Series representation for arbitrary positive order}\label{sec:5_2}

\cfclear
\begin{athm}{lemma}{lem:ver2_pre}
Let $s \in \N$. Then it holds for all $k\in\Z$ that
\begin{equation}\label{eq:5.13}
\lim_{z \to s } \Kern_z(k) = 
\frac{ \1_{\{1,2,\ldots,s\}}(k) \, (-1)^{k+1} \gamfn(2s +1) }{ \gamfn(1+s+k) \gamfn(1+s-k) }
\end{equation}
\cfout.
\end{athm}

\begin{aproof}
First, note that \cref{prop:gamma} ensures that for all $k \in \Z$ it holds that
\begin{equation}
\lim_{z\to s} \Kern_z(k) = \lim_{z\to s} \frac{ \1_{\Z\backslash\{0\}}(k) \, (-1)^{k+1} \gamfn(2s+1) }{ \gamfn( 1 + s + k ) \gamfn( 1 + s - k ) }
\end{equation}
\cfload.
Combining this and the fact that \cref{def:gamma} ensures that for all $z \in \{\ldots,-2,-1,0\}$ it holds that $\nicefrac{1}{\gamfn(z)} = 0$ establishes \cref{eq:5.13}.
\end{aproof}

\cfclear
\begin{athm}{lemma}{lem:ver2}
It holds for all $s \in \N$, $u \in \ellTwo$, $n \in \Z$ that
\begin{equation}\label{this0}
\pr[\big]{ \fracLap u }(n) = \sum_{k=0}^{2s} (-1)^{k-s} \binom{2s}{k} u(n-s+k) = \lim_{z \to s} \br[\Bigg]{  \sum_{k\in\Z} \Kern_z(n-k) \pr[\big]{ u(n) - u(k) } } \ifnocf.
\end{equation}
\cfload.
\end{athm}

\begin{aproof}
First, note thatthe fact that for all $a \in \N$, $b \in \{0,1,\ldots,a\}$ it holds that $\binom{a}{b} = \binom{a}{a-b}$ assures that for all $s\in\N$, $u\in\ellTwo$, $n\in\Z$ it holds that
\begin{align}\label{this1}
& \sum_{k=0}^{2s} (-1)^{k-s} \binom{2s}{k} u(n-s+k)
= \sum_{k=-s}^s (-1)^k \binom{2s}{s+k} u(n-k) \nonumber \\
& \quad = \sum_{k=-s}^s \br[\bigg]{ \frac{ (-1)^k \gamfn( 2s+1 ) }{ \gamfn( 1+s+k ) \gamfn(1+s-k) } } u(n-k) \\
& \quad = \br[\bigg]{ \frac{ \gamfn(2s+1) }{ \gamfn(1+s) \gamfn(1+s) } } u(n) -  \sum_{k=1}^s \br[\bigg]{ \frac{ (-1)^{k+1} \gamfn( 2s+1 ) }{ \gamfn( 1+s+k ) \gamfn(1+s-k) } } \pr[\big]{ u(n-k) + u(n+k) }  \nonumber
\end{align}
\cfload.
Next, observe that \cref{lem:sum} ensures that for all $z \in (0,\infty) \backslash \N$, $u \in \ellTwo$, $n\in\Z$ it holds that
\begin{align}
\sum_{k\in\Z} \Kern_z(n-k) \pr[\big]{ u(n) - u(k) } 
& = u(n) \sum_{k\in\Z} \Kern_z(n-k) - \sum_{k\in\Z} \Kern_z(n-k) u(k) \nonumber \\
& = \br[\bigg]{ \frac{4^z \gamfn(\nicefrac{1}{2} + z)}{\sqrt{\pi} \gamfn(1+z) } } u(n) - \sum_{k\in\Z} \Kern_z(n-k) u(k) \ifnocf.
\end{align}
This, \cref{lem:kern_bd_i2,lem:kern_bd_i1} of \cref{lem:kern_bd}, \cref{lem:ver2_pre}, and, e.g., Rudin \cite[Theorem 7.17]{rudin1964principles} guarantee that for all $s\in\N$, $u \in \ellTwo$, $n\in\Z$ it holds that
\begin{align}\label{this2}
& \lim_{z \to s} \br[\Bigg]{  \sum_{k\in\Z} \Kern_s(n-k) \pr[\big]{ u(n) - u(k) } }
= \br[\bigg]{ \lim_{z\to s} \frac{4^z \gamfn(\nicefrac{1}{2} + z)}{\sqrt{\pi} \gamfn(1+z) } } u(n) - \lim_{z\to s} \left[\sum_{k\in\Z} \Kern_z(n-k) u(k)\right] \nonumber \\
& \quad = \br[\bigg]{ \frac{4^s \gamfn(\nicefrac{1}{2} + s)}{\sqrt{\pi} \gamfn(1+s) } } u(n) - \lim_{z\to s} \left[ \sum_{k\in\N} \Kern_z(k) \pr[\big]{ u(n-k) + u(n+k) } \right] \nonumber \\
& \quad = \br[\bigg]{ \frac{4^s \gamfn(\nicefrac{1}{2} + s)}{\sqrt{\pi} \gamfn(1+s) } } u(n) - \sum_{k\in\N} \br[\Big]{ \lim_{z\to s}\Kern_z(k) } \pr[\big]{ u(n-k) + u(n+k) } \\
& \quad = \br[\bigg]{ \frac{4^s \gamfn(\nicefrac{1}{2} + s)}{\sqrt{\pi} \gamfn(1+s) } } u(n) - \sum_{k\in\N} \br[\bigg]{ \frac{ \1_{\{1,2,\ldots,s\}}(k) \, (-1)^{k+1} \gamfn(2s+1)}{\gamfn(1+s+k)\gamfn(1+s-k)} } \pr[\big]{ u(n-k) + u(n+k) } \nonumber \\
& \quad = \br[\bigg]{ \frac{4^s \gamfn(\nicefrac{1}{2} + s)}{\sqrt{\pi} \gamfn(1+s) } } u(n) - \sum_{k=1}^s \br[\bigg]{ \frac{ (-1)^{k+1} \gamfn(2s+1)}{\gamfn(1+s+k)\gamfn(1+s-k)} } \pr[\big]{ u(n-k) + u(n+k) }. \nonumber 
\end{align}
In addition, note that the fact that \cref{def:gamma} implies that for all $z\in\C$ with $\real(z) \in \R \backslash\{\ldots,-2,\allowbreak -1, \allowbreak 0\}$ it holds that $z\gamfn(z) = \gamfn(1+z)$ and the Legendre duplication formula 
(cf., e.g., \cite[Eq.\ (6.1.18), Page 256]{abramowitz1970handbook}) demonstrate that for all $s\in\N$ it holds that
\begin{align}
\frac{4^s \gamfn(\nicefrac{1}{2} + s)}{\sqrt{\pi} \gamfn(1+s) }
& = \frac{4^s \gamfn(\nicefrac{1}{2} + s)}{\sqrt{\pi} \gamfn(1+s) } \cdot \frac{s\gamfn(s)}{s\gamfn(s)}
= \frac{4^s \br[\big]{ 2^{1-2s} \sqrt{\pi} s \gamfn(2s) }}{\sqrt{\pi} \gamfn(s) \gamfn(1+s) }
= \frac{ \gamfn(2s+1) }{ \gamfn(1+s) \gamfn(1+s) }.
\end{align}
Combining this, \cref{this1}, \cref{this2}, and \cref{lem:ver1} proves \cref{this0}.
\end{aproof}

\cfclear
\begin{athm}{theorem}{cor:final}
It holds for all $s\in(0,\infty)$, $u \in \ellTwo$, $n\in\Z$ that
\begin{equation}\label{eq:frac_rep_final}
\pr[\big]{ \fracLap u }(n) = \lim_{z \to s} \br[\Bigg]{  \sum_{k\in\Z} \Kern_z(n-k) \pr[\big]{ u(n) - u(k) } } \ifnocf.
\end{equation}
\cfload.
\end{athm}

\begin{aproof}
Note that combining \cref{lem:rep1,lem:ver2} establishes \cref{eq:frac_rep_final}. 
\end{aproof}


\section{Conclusions and future endeavors}\label{sec:6}

\subsection{Concluding remarks}

In this article we developed novel results regarding real-valued positive fractional powers of the discrete Laplace operator.
In particular, we defined a discrete fractional Laplace operator for arbitrary real-valued positive powers (cf.\ \cref{def:frac_lap}) 
and then developed its series representation (cf.\ \cref{cor:final}).
This latter task was primarily accomplished through the development of two sets of results.
First, we constructed the series representation for positive integer powers of the discrete Laplace operator (cf.\ \cref{lem:ver1,lem:ver2}).
Next, we developed series representations for positive non-integer powers of the discrete Laplace operator (cf.\ \cref{lem:rep1}). 
The main result of the article (cf.\ \cref{cor:final}) is obtained by showing that the series representations 
obtained in each of the previous steps in fact coincide.

The main results developed---i.e., the results of \cref{sec:5}---required numerous preliminary results from various areas of mathematics. 
The results in \cref{sec:2} are of a functional analysis flavor and allow for a beautiful description of 
important properties of strongly continuous semigroups.
These results were combined with results from discrete harmonic analysis in \cref{sec:3} in order to define 
and study the discrete fractional Laplace operator.
Since the presented definition of this operator (cf.\ \cref{def:frac_lap}) employs a so-called semigroup 
language, it was imperative that all novel mathematical objects are determined to be well 
defined in $\ellTwo$ (cf.\ \cref{def:l2}).
Finally, \cref{sec:4} provides a detailed study of the proposed fractional kernel function (cf.\ \cref{def:kern}) which is necessary for the development of the coefficients of the 
series representations presented in \cref{sec:5}. 
Therein, it is shown that the proposed fractional kernel function is well-defined, 
symmetric, and continuous for all $s \in (0,\infty) \backslash \N$.
It is later shown in \cref{lem:ver2_pre} that the values $s\in\N$ are in fact removable singularities.

As a final remark, we wish to emphasize the importance of the presented results (for a specific physical motivation, see \cref{sec:motivate} above). 
Due to the rapidly growing interest in problems related to fractional calculus, there is a need to 
determine the validity of including fractional operators into existing models.
The study of the discrete fractional Laplace operator, or its continuous counterpart, for the case 
when $s\in(0,1)$ is well-understood and often used in physical sciences.
In this setting, the operator may be used to model \textit{superdiffusive} phenomena \cite{padgett2020anomalous}.
Moreover, there have been rigorous studies of the operator in this parameter regime which 
demonstrate that such considerations are well-defined and well-behaved.
As such, it is natural to attempt to extend these studies to the case when $s\in(1,\infty)$, as well.
The current article has demonstrated that such extensions are indeed well-defined in the discrete case.
In addition, \cref{cor:final} shows that one may potentially use the existing understanding of 
the case of positive integer powers of the discrete Laplace operator to provide some much needed 
physical intuition to the discrete fractional Laplace operator.
However, there are still numerous unanswered questions regarding important properties of 
these operators and we will outline a few such open problems and research directions in \cref{sec:6_2} below.

\subsection{Related future endeavors}\label{sec:6_2}

First and foremost, there is a need to continue the analytical work presented herein in order to 
obtain a better understanding of the discrete fractional Laplace operator.
In this article, we have considered the setting where all objects are defined in $\ellTwo$;
however, this is not always the appropriate setting for physically relevant problems.
As such, we intend to extend our study to the situation where the underlying function spaces have less 
regularity (e.g., discrete H{\"o}lder spaces) and develop standard regularity estimates.
We also intend to develop similar series representations for the case of real-valued negative exponents.
Such representations are highly important for studying fractional Poisson-like problems, 
as they provide representations of the solutions to these problems.
Finally, we hope to develop an understanding of the spectral properties of the discrete fractional Laplace operator.
While it is well-known that the discrete Laplace operator (cf.\ \cref{def:lap}) has purely absolutely 
continuous spectra (cf., e.g., Dutkay and Jorgensen \cite{dutkay2010spectral}), to the authors' 
knowledge this has not been rigorously proven in the case of the discrete fractional Laplace operator (although the result is expected).
Demonstrating this to be the case is of utmost importance and will have far-reaching implications in 
mathematics and physics (for clarification, see the techniques outlined in Liaw \cite{liaw2013approach}).

The proposed discrete fractional Laplace operator is also of interest due to its importance in
multi-physical sciences.
An example of an avenue for future research is transport in turbulent plasmas. It has been experimentally 
observed that heat and particle transport in turbulent plasmas 
is \textit{non-local} (i.e., \textit{anomalous}) in nature 
(cf., e.g., \cite{shalchi2011magnetic,del2005nondiffusive,%
gentle1995experimental}). 
Comparison between transport models using the fractional Laplace operator and experimental results 
have demonstrated that electron transport in turbulent fusion plasmas is characterized by fractional 
exponents in the range 
$s\in(0.6,1)$, which indicates superdiffusive behavior 
\cite{kullberg2014comparison,kullberg2013isotropic}. 
Moreover, a generalized approach to modeling anomalous diffusive transport in turbulent plasmas 
employs diffusion-type equations where fractional derivatives occur in both space and time 
(cf., e.g., del Castillo-Negrete et al.\ \cite{del2004fractional}). 
Using the series representations presented herein, we intend to show that the fractional derivative 
in time can be incorporated into the spatial derivative, which can greatly simplify such equations. 
Finally, better modeling explorations of subdiffusive and superdiffusive systems and 
fractional quenching-combustion phenomena will also be within our continuing endeavors (cf., e.g., \cite{MR4066721,padgett2018quenching}).


\section*{Acknowledgments}

The second author acknowledges funding by the National Science Foundation (NSF 1903450) and the Department of Energy (DE-SC0021284). 
The third author acknowledges funding by the National Science Foundation (NSF 1903450).
The fourth author would like to thank the College of Arts and Sciences at Baylor University for partial support 
through a research leave award.

\bibliographystyle{acm}
\bibliography{bibfile}

\begin{thebibliography}{10}

\bibitem{abramowitz1970handbook}
{\sc Abramowitz, M., and Stegun, I.~A.}, Eds.
\newblock {\em Handbook of mathematical functions with formulas, graphs, and
  mathematical tables}.
\newblock Dover Publications, Inc., New York, 1992.
\newblock Reprint of the 1972 edition.

\bibitem{aizenman1993localization}
{\sc Aizenman, M., and Molchanov, S.}
\newblock Localization at large disorder and at extreme energies: {A}n
  elementary derivations.
\newblock {\em Communications in Mathematical Physics 157}, 2 (1993), 245--278.

\bibitem{MR2215610}
{\sc Aizenman, M., Sims, R., and Warzel, S.}
\newblock Absolutely continuous spectra of quantum tree graphs with weak
  disorder.
\newblock {\em Comm. Math. Phys. 264}, 2 (2006), 371--389.

\bibitem{PhysRev.109.1492}
{\sc Anderson, P.~W.}
\newblock Absence of diffusion in certain random lattices.
\newblock {\em Phys. Rev. 109\/} (Mar 1958), 1492--1505.

\bibitem{arendt2018fractional}
{\sc Arendt, W., ter Elst, A. F.~M., and Warma, M.}
\newblock Fractional powers of sectorial operators via the
  {D}irichlet-to-{N}eumann operator.
\newblock {\em Comm. Partial Differential Equations 43}, 1 (2018), 1--24.

\bibitem{artin2015gamma}
{\sc Artin, E.}
\newblock {\em The gamma function}.
\newblock Translated by Michael Butler. Athena Series: Selected Topics in
  Mathematics. Holt, Rinehart and Winston, New York-Toronto-London, 1964.

\bibitem{brandes2003anderson}
{\sc Brandes, T., and Kettemann, S.}
\newblock {\em {Anderson localization and its ramifications: Disorder, phase
  coherence, and electron correlations}}, vol.~630.
\newblock Springer Science \& Business Media, 2003.

\bibitem{bucur2016nonlocal}
{\sc Bucur, C., and Valdinoci, E.}
\newblock {\em Nonlocal diffusion and applications}, vol.~20 of {\em Lecture
  Notes of the Unione Matematica Italiana}.
\newblock Springer, [Cham]; Unione Matematica Italiana, Bologna, 2016.

\bibitem{caffarelli2007extension}
{\sc Caffarelli, L., and Silvestre, L.}
\newblock An extension problem related to the fractional {L}aplacian.
\newblock {\em Comm. Partial Differential Equations 32}, 7-9 (2007),
  1245--1260.

\bibitem{cartan2017differential}
{\sc Cartan, H., Moore, J., Husemoller, D., and Maestro, K.}
\newblock {\em Differential Calculus on Normed Spaces: A Course in Analysis}.
\newblock CreateSpace Independent Publishing Platform, 2017.

\bibitem{chen2018extension}
{\sc Chen, Y.~K., Lei, Z., and Wei, C.~H.}
\newblock Extension problems related to the higher order fractional
  {L}aplacian.
\newblock {\em Acta Math. Sin. (Engl. Ser.) 34}, 4 (2018), 655--661.

\bibitem{ciaurri2017harmonic}
{\sc Ciaurri, O., Gillespie, T.~A., Roncal, L., Torrea, J.~L., and Varona,
  J.~L.}
\newblock Harmonic analysis associated with a discrete {L}aplacian.
\newblock {\em J. Anal. Math. 132\/} (2017), 109--131.

\bibitem{ciaurri2018nonlocal}
{\sc Ciaurri, O., Roncal, L., Stinga, P.~R., Torrea, J.~L., and Varona, J.~L.}
\newblock Nonlocal discrete diffusion equations and the fractional discrete
  {L}aplacian, regularity and applications.
\newblock {\em Adv. Math. 330\/} (2018), 688--738.

\bibitem{del2004fractional}
{\sc del Castillo-Negrete, D., Carreras, B., and Lynch, V.}
\newblock Fractional diffusion in plasma turbulence.
\newblock {\em Physics of Plasmas 11}, 8 (2004), 3854--3864.

\bibitem{del2005nondiffusive}
{\sc del Castillo-Negrete, D., Carreras, B., and Lynch, V.}
\newblock Nondiffusive transport in plasma turbulence: a fractional diffusion
  approach.
\newblock {\em Physical review letters 94}, 6 (2005), 065003.

\bibitem{dragomir2004semi}
{\sc Dragomir, S.~S.}
\newblock {\em Semi-inner products and applications}.
\newblock Nova Science Publishers, Inc., Hauppauge, NY, 2004.

\bibitem{driverfunctional}
{\sc Driver, B.~K.}
\newblock Functional analysis tools with examples.

\bibitem{dutkay2010spectral}
{\sc Dutkay, D.~E., and Jorgensen, P. E.~T.}
\newblock Spectral theory for discrete {L}aplacians.
\newblock {\em Complex Anal. Oper. Theory 4}, 1 (2010), 1--38.

\bibitem{felli2018unique}
{\sc Felli, V., and Ferrero, A.}
\newblock Unique continuation principles for a higher order fractional
  {L}aplace equation.
\newblock {\em Nonlinearity 33}, 8 (2020), 4133--4190.

\bibitem{frohlich1983absence}
{\sc Fr{\"o}hlich, J., and Spencer, T.}
\newblock {Absence of diffusion in the Anderson tight binding model for large
  disorder or low energy}.
\newblock {\em Communications in Mathematical Physics 88}, 2 (1983), 151--184.

\bibitem{gale2013extension}
{\sc Gal\'{e}, J.~E., Miana, P.~J., and Stinga, P.~R.}
\newblock Extension problem and fractional operators: semigroups and wave
  equations.
\newblock {\em J. Evol. Equ. 13}, 2 (2013), 343--368.

\bibitem{garcia2019strong}
{\sc Garc\'{\i}a-Ferrero, M.~A., and R\"{u}land, A.}
\newblock Strong unique continuation for the higher order fractional
  {L}aplacian.
\newblock {\em Math. Eng. 1}, 4 (2019), 715--774.

\bibitem{gentle1995experimental}
{\sc Gentle, K., Bravenec, R., Cima, G., Gasquet, H., Hallock, G., Phillips,
  P., Ross, D., Rowan, W., Wootton, A., Crowley, T., et~al.}
\newblock An experimental counter-example to the local transport paradigm.
\newblock {\em Physics of Plasmas 2}, 6 (1995), 2292--2298.

\bibitem{giles1967classes}
{\sc Giles, J.~R.}
\newblock Classes of semi-inner-product spaces.
\newblock {\em Trans. Amer. Math. Soc. 129\/} (1967), 436--446.

\bibitem{PhysRevLett.102.085002}
{\sc Hou, L.-J., Piel, A., and Shukla, P.~K.}
\newblock Self-diffusion in 2{D} dusty-plasma liquids: {N}umerical-simulation
  results.
\newblock {\em Phys. Rev. Lett. 102\/} (Feb 2009), 085002.

\bibitem{MR1779620}
{\sc Jak\v{s}i\'{c}, V., and Last, Y.}
\newblock Spectral structure of {A}nderson type {H}amiltonians.
\newblock {\em Invent. Math. 141}, 3 (2000), 561--577.

\bibitem{MR2219273}
{\sc Jak\v{s}i\'{c}, V., and Last, Y.}
\newblock Simplicity of singular spectrum in {A}nderson-type {H}amiltonians.
\newblock {\em Duke Math. J. 133}, 1 (2006), 185--204.

\bibitem{jones1}
{\sc Jones, T.~F., Padgett, J.~L., and Sheng, Q.}
\newblock Intrinsic properties of strongly continuous fractional semigroups in
  normed vector spaces.
\newblock {\em arXiv preprint arXiv:2012.11092\/} (2020).

\bibitem{kelley2001difference}
{\sc Kelley, W.~G., and Peterson, A.~C.}
\newblock {\em Difference equations}, second~ed.
\newblock Harcourt/Academic Press, San Diego, CA, 2001.
\newblock An introduction with applications.

\bibitem{kostadinova2021active}
{\sc Kostadinova, E.~G., Banka, R., Padgett, J.~L., Liaw, C.~D., Matthews,
  L.~S., and Hyde, T.~W.}
\newblock Active turbulence in a dusty plasma monolayer.
\newblock {\em arXiv preprint arXiv:2102.09344\/} (2021).

\bibitem{kostadinova2017delocalization}
{\sc Kostadinova, E.~G., Busse, K., Ellis, N., Padgett, J.~L., Liaw, C.~D.,
  Matthews, L.~S., and Hyde, T.~W.}
\newblock Delocalization in infinite disordered two-dimensional lattices of
  different geometry.
\newblock {\em Physical Review B 96}, 23 (2017), 235408.

\bibitem{PhysRevResearch.2.043375}
{\sc Kostadinova, E.~G., Padgett, J.~L., Liaw, C.~D., Matthews, L.~S., and
  Hyde, T.~W.}
\newblock Numerical study of anomalous diffusion of light in semicrystalline
  polymer structures.
\newblock {\em Phys. Rev. Research 2\/} (Dec 2020), 043375.

\bibitem{kramer1993localization}
{\sc Kramer, B., and MacKinnon, A.}
\newblock Localization: theory and experiment.
\newblock {\em Reports on Progress in Physics 56}, 12 (1993), 1469.

\bibitem{kreyszig1978introductory}
{\sc Kreyszig, E.}
\newblock {\em Introductory functional analysis with applications}.
\newblock Wiley Classics Library. John Wiley \& Sons, Inc., New York, 1989.

\bibitem{kullberg2013isotropic}
{\sc Kullberg, A., del Castillo-Negrete, D., Morales, G., and Maggs, J.}
\newblock Isotropic model of fractional transport in two-dimensional bounded
  domains.
\newblock {\em Physical Review E 87}, 5 (2013), 052115.

\bibitem{kullberg2014comparison}
{\sc Kullberg, A., Morales, G., and Maggs, J.}
\newblock Comparison of a radial fractional transport model with tokamak
  experiments.
\newblock {\em Physics of Plasmas 21}, 3 (2014), 032310.

\bibitem{lagendijk2009fifty}
{\sc Lagendijk, A., Van~Tiggelen, B., and Wiersma, D.~S.}
\newblock Fifty years of {A}nderson localization.
\newblock {\em Phys. Today 62}, 8 (2009), 24--29.

\bibitem{poin}
{\sc Lebedev, V.~I., and Agoshkov, V.~I.}
\newblock {\em {Poincar{\'e}-Steklov operators and their applications in
  analysis}}.
\newblock Akad. Nauk SSSR, Vychisl. Tsentr, Moscow, 1983.

\bibitem{liaw2013approach}
{\sc Liaw, C.}
\newblock Approach to the extended states conjecture.
\newblock {\em Journal of Statistical Physics 153}, 6 (2013), 1022--1038.

\bibitem{lischke2020fractional}
{\sc Lischke, A., Pang, G., Gulian, M., and et~al.}
\newblock What is the fractional {L}aplacian? {A} comparative review with new
  results.
\newblock {\em J. Comput. Phys. 404\/} (2020), 109009, 62.

\bibitem{PhysRevE.75.016405}
{\sc Liu, B., and Goree, J.}
\newblock {Superdiffusion in two-dimensional Yukawa liquids}.
\newblock {\em Phys. Rev. E 75\/} (Jan 2007), 016405.

\bibitem{PhysRevLett.100.055003}
{\sc Liu, B., and Goree, J.}
\newblock {Superdiffusion and non-Gaussian statistics in a driven-dissipative
  2D dusty plasma}.
\newblock {\em Phys. Rev. Lett. 100\/} (Feb 2008), 055003.

\bibitem{lumer1961semi}
{\sc Lumer, G.}
\newblock Semi-inner-product spaces.
\newblock {\em Trans. Amer. Math. Soc. 100\/} (1961), 29--43.

\bibitem{meerschaert2011stochastic}
{\sc Meerschaert, M.~M., and Sikorskii, A.}
\newblock {\em Stochastic models for fractional calculus}, second~ed., vol.~43
  of {\em De Gruyter Studies in Mathematics}.
\newblock De Gruyter, Berlin, 2019.

\bibitem{meichsner2017fractional}
{\sc Meichsner, J., and Seifert, C.}
\newblock {Fractional powers of non-negative operators in Banach spaces via the
  Dirichlet-to-Neumann operator}.
\newblock {\em arXiv preprint arXiv:1704.01876\/} (2017).

\bibitem{meichsner2020harmonic}
{\sc Meichsner, J., and Seifert, C.}
\newblock On the harmonic extension approach to fractional powers in {B}anach
  spaces.
\newblock {\em Fract. Calc. Appl. Anal. 23}, 4 (2020), 1054--1089.

\bibitem{PhysRevLett.96.015003}
{\sc Nunomura, S., Samsonov, D., Zhdanov, S., and Morfill, G.}
\newblock Self-diffusion in a liquid complex plasma.
\newblock {\em Phys. Rev. Lett. 96\/} (Jan 2006), 015003.

\bibitem{olver2010nist}
{\sc Olver, F. W.~J., Lozier, D.~W., Boisvert, R.~F., and Clark, C.~W.}, Eds.
\newblock {\em N{IST} handbook of mathematical functions}.
\newblock U.S. Department of Commerce, National Institute of Standards and
  Technology, Washington, DC; Cambridge University Press, Cambridge, 2010.
\newblock With 1 CD-ROM (Windows, Macintosh and UNIX).

\bibitem{https://doi.org/10.1002/ctpp.200910089}
{\sc Ott, T., and Bonitz, M.}
\newblock Anomalous and {F}ickian diffusion in two-dimensional dusty plasmas.
\newblock {\em Contributions to Plasma Physics 49}, 10 (2009), 760--764.

\bibitem{PhysRevE.78.026409}
{\sc Ott, T., Bonitz, M., Donk\'o, Z., and Hartmann, P.}
\newblock {Superdiffusion in quasi-two-dimensional Yukawa liquids}.
\newblock {\em Phys. Rev. E 78\/} (Aug 2008), 026409.

\bibitem{padgett2018quenching}
{\sc Padgett, J.~L.}
\newblock The quenching of solutions to time-space fractional {K}awarada
  problems.
\newblock {\em Comput. Math. Appl. 76}, 7 (2018), 1583--1592.

\bibitem{padgett2020analysis}
{\sc Padgett, J.~L.}
\newblock Analysis of an approximation to a fractional extension problem.
\newblock {\em BIT 60}, 3 (2020), 715--739.

\bibitem{padgett2020anomalous}
{\sc Padgett, J.~L., Kostadinova, E.~G., Liaw, C.~D., Busse, K., Matthews,
  L.~S., and Hyde, T.~W.}
\newblock Anomalous diffusion in one-dimensional disordered systems: a discrete
  fractional {L}aplacian method.
\newblock {\em J. Phys. A 53}, 13 (2020), 135205, 21.

\bibitem{pastur1980spectral}
{\sc Pastur, L.~A.}
\newblock Spectral properties of disordered systems in the one-body
  approximation.
\newblock {\em Communications in Mathematical Physics 75}, 2 (1980), 179--196.

\bibitem{pozrikidis2018fractional}
{\sc Pozrikidis, C.}
\newblock {\em The fractional {L}aplacian}.
\newblock CRC Press, Boca Raton, FL, 2016.

\bibitem{ros2014local}
{\sc Ros-Oton, X., and Serra, J.}
\newblock Local integration by parts and {P}ohozaev identities for higher order
  fractional {L}aplacians.
\newblock {\em Discrete Contin. Dyn. Syst. 35}, 5 (2015), 2131--2150.

\bibitem{rudin1964principles}
{\sc Rudin, W.}
\newblock {\em Principles of mathematical analysis}, third~ed.
\newblock McGraw-Hill Book Co., New York-Auckland-D\"{u}sseldorf, 1976.
\newblock International Series in Pure and Applied Mathematics.

\bibitem{shalchi2011magnetic}
{\sc Shalchi, A.}
\newblock Magnetic field line random walk in two-dimensional turbulence:
  {M}arkovian diffusion versus superdiffusion.
\newblock {\em Contributions to Plasma Physics 51}, 10 (2011), 920--930.

\bibitem{soderlind2006logarithmic}
{\sc S\"{o}derlind, G.}
\newblock The logarithmic norm. {H}istory and modern theory.
\newblock {\em BIT 46}, 3 (2006), 631--652.

\bibitem{PhysRevX.6.021021}
{\sc Strickler, T.~S., Langin, T.~K., McQuillen, P., Daligault, J., and
  Killian, T.~C.}
\newblock Experimental measurement of self-diffusion in a strongly coupled
  plasma.
\newblock {\em Phys. Rev. X 6\/} (May 2016), 021021.

\bibitem{tricomi1951asymptotic}
{\sc Tricomi, F.~G., and Erd\'{e}lyi, A.}
\newblock The asymptotic expansion of a ratio of gamma functions.
\newblock {\em Pacific J. Math. 1\/} (1951), 133--142.

\bibitem{doi:10.1063/1.1449888}
{\sc Vaulina, O.~S., and Vladimirov, S.~V.}
\newblock Diffusion and dynamics of macro-particles in a complex plasma.
\newblock {\em Physics of Plasmas 9}, 3 (2002), 835--840.

\bibitem{vazquez2017mathematical}
{\sc V\'{a}zquez, J.~L.}
\newblock The mathematical theories of diffusion: nonlinear and fractional
  diffusion.
\newblock In {\em Nonlocal and nonlinear diffusions and interactions: new
  methods and directions}, vol.~2186 of {\em Lecture Notes in Math.} Springer,
  Cham, 2017, pp.~205--278.

\bibitem{yang2013higher}
{\sc Yang, R.}
\newblock On higher order extensions for the fractional {L}aplacian.
\newblock {\em arXiv preprint arXiv:1302.4413\/} (2013).

\bibitem{MR4066721}
{\sc Zhu, L., and Sheng, Q.}
\newblock A note on the adaptive numerical solution of a {R}iemann-{L}iouville
  space-fractional {K}awarada problem.
\newblock {\em J. Comput. Appl. Math. 374\/} (2020), 112714, 14pp.

\end{thebibliography}

\end{document}